\documentclass[12pt, leqno]{amsart}
\usepackage{graphicx}
\usepackage{amssymb,amsmath,amsthm,amscd}
\usepackage{mathrsfs}
\usepackage{enumerate}
\usepackage{enumitem}
\usepackage{verbatim}
\usepackage[usenames,dvipsnames]{color}
\usepackage[colorlinks=true, pdfstartview=FitV,
linkcolor=blue,citecolor=blue,urlcolor=blue]{hyperref}
\usepackage[all]{xy}
\usepackage{tikz}
\usepackage{tikz-cd}
\usetikzlibrary{matrix, arrows}
\usepackage{stmaryrd} 
\usepackage{dsfont}
\usepackage{bm} 
\usepackage{scalerel}[2016/12/29] 
\usepackage{amsthm}
\usepackage{adjustbox}
\usepackage{float}

\allowdisplaybreaks[4]

\DeclareMathSymbol{\shortminus}{\mathbin}{AMSa}{"39}

\newcommand{\nc}{\newcommand}

\numberwithin{equation}{section}

\usepackage[colorinlistoftodos]{todonotes}
\usepackage{mathtools}
\mathtoolsset{showonlyrefs,showmanualtags}
\usepackage{empheq}
\theoremstyle{plain}
\newtheorem{lem}{Lemma}[section]
\newtheorem{pro}[lem]{Proposition}
\newtheorem{thm}[lem]{Theorem}
\newtheorem{cor}[lem]{Corollary}
\newtheorem{defi}[lem]{Definition}
\newtheorem{conjecture}[lem]{Conjecture}

\newcommand{\Pro}{\begin{pro}}
	\newcommand{\enpro}{\end{pro}}
\newcommand{\Lem}{\begin{lem}}
	\newcommand{\enlem}{\end{lem}}
\newcommand{\Thm}{\begin{thm}}
	\newcommand{\enthm}{\end{thm}}
\newcommand{\Cor}{\begin{cor}}
	\newcommand{\encor}{\end{cor}}
\newcommand{\Defi}{\begin{defi}}
	\newcommand{\enDefi}{\end{defi}}
\newcommand{\Proof}{\begin{proof}}
	\newcommand{\enproof}{\end{proof}}

\theoremstyle{definition} 

\newtheorem{Convention}[lem]{Convention}
\newtheorem{notat}[lem]{Notation}

\theoremstyle{remark}
\newtheorem{rem}[lem]{Remark}

\newcommand{\Conv}{\begin{Convention}}
	\newcommand{\enconv}{\end{Convention}}
\nc{\Rem}{\begin{rem}}
	\nc{\enrem}{\end{rem}}

\newcommand{\arxiv}[1]{\href{http://arxiv.org/abs/#1}{\tt arXiv:\nolinkurl{#1}}}

\nc{\rmkend}{\hfill$\triangledown$}
\nc{\defend}{\hfill$\triangle$}
\nc{\Tr}{\on{Tr}}
\nc{\id}{\on{id}}
\nc{\mi}{^{-1}}
\nc{\indx}{\mathbb{I}}
\nc{\Sym}{\on{Sym}}

\nc{\fext}[2]{{#1}[\negthinspace[#2]\negthinspace]}
\nc{\gext}[2]{{#1}(\negthinspace(#2)\negthinspace)}

\newcommand{\on}{\operatorname}
\nc{\be}{\begin{enumerate}}
	\nc{\ee}{\end{enumerate}}
\newcommand{\eq}{\begin{equation}}
	\newcommand{\eneq}{\end{equation}}
\nc{\bc}{\begin{cases}}
	\nc{\ec}{\end{cases}}
\newcommand{\eqn}{\begin{eqnarray*}}
	\newcommand{\eneqn}{\end{eqnarray*}}
\newcommand{\ba}{\begin{array}}
	\newcommand{\ea}{\end{array}}

\newcommand{\commentout}[1]{}

\newcommand{\C}{{\mathbb C}}
\newcommand{\Q}{\mathbb {Q}}
\newcommand{\Z}{{\mathbb Z}}
\newcommand{\N}{{\mathbb N}}

\nc{\Ad}{\operatorname{Ad}}
\nc{\gr}{\on{gr}}

\nc{\Aut}{\operatorname{Aut}}

\nc{\coker}{\operatorname{coker}}

\nc{\Img}{\on{Im}}
\nc{\res}{\on{res}} 

\nc{\bl}{\bigl(}
\nc{\br}{\bigr)}

\nc{\Ya}{\widetilde{Y}_{\hbar}}
\nc{\Yu}{\mathbf{Y}_{\mu}}
\nc{\Yinf}{\mathbf{Y}_{\infty}}
\nc{\hcor}{\C[\hbar]}
\nc{\Tya}{{}^{\mathsf{tw}}\widetilde{Y}_{\mu,\hbar}} 
\nc{\Tyazero}{{}^{\mathsf{tw}}\widetilde{Y}_{0,\hbar}} 
\nc{\Ty}{{}^{\mathsf{tw}}\widetilde{{Y}}_{\hbar}} 
\nc{\sTy}{{}^{\mathsf{tw}}{{Y}}_{\hbar}} 
\nc{\dya}{{}^{\mathsf{tw}}{\mathbf{Y}}_{{\mu}}} 
\nc{\dyatr}{{}^{\mathsf{tw}}{\mathbf{Y}}_{{\mu}}^{{\lambda}}} 
\nc{\dyatrzero}{{}^{\mathsf{tw}}{\mathbf{Y}}_{{0}}^{{\lambda}}} 
\nc{\dyazero}{{}^{\mathsf{tw}}{\mathbf{Y}}_0}
\nc{\yazero}{{\mathbf{Y}}_{0}}
\nc{\ttya}{{}^{\mathsf{tw}}\widetilde{{Y}}_{\mu,\hbar}^\lambda} 
\nc{\sttya}{{}^{\mathsf{tw}}{{Y}}_{{\mu},\hbar}^{{\lambda}}}
\nc{\Tyar}{{}^{\mathsf{tw}}\widetilde{{Y}}_{\mu}}
\nc{\sTya}{{}^{\mathsf{tw}}{Y}_{{\mu},\hbar}} 
\nc{\sTyar}{{}^{\mathsf{tw}}{{Y}}_{{\mu}}}
\nc{\Rty}{{}^{\mathsf{tw}}\widetilde{{\mathbf{Y}}}_{\mu}}

\nc{\tYrtt}{{}^{\mathsf{tw}}\widetilde{Y}_{\hbar}^{\mathsf{rtt}}}
\nc{\Yartt}{\widetilde{Y}_{\hbar}^{\mathsf{rtt}}}
\nc{\sYartt}{{Y}_{\hbar}^{\mathsf{rtt}}} 
\nc{\bYartt}{\widetilde{\mathbf{Y}}_{}^{\mathsf{rtt}}}
\nc{\bsYartt}{\mathbf{Y}_{}^{\mathsf{rtt}}}

\nc{\h}{\hbar}
\nc{\diag}{\on{diag}}
\nc{\lat}{\Lambda^\vee}
\nc{\dlat}{\Lambda}
\nc{\splat}{\Lambda^{2+}}
\nc{\slat}{\overline{\Lambda}}
\nc{\uh}{\underline{h}}
\nc{\uz}{\underline{z}}
\nc{\eindx}{\tilde{\indx}}
\nc{\ai}{\alpha_i^\vee}
\nc{\ei}{\epsilon_i^\vee}
\nc{\gln}{\mathfrak{gl}_n}
\nc{\sln}{\mathfrak{sl}_n}
\nc{\tih}{\widetilde{h}}
\newcommand{\half}{{\textstyle \frac{1}{2}}}
\nc{\diff}{\mathbf{D}_\lambda}
\nc{\hdiff}{{D}_{\mu,\hbar}^\lambda}    
\nc{\hdiffs}{\breve{D}_{\mu,\hbar}^\lambda} 
\nc{\hdiffzero}{{D}_{0,\hbar}^\lambda}
\nc{\LHS}{\on{LHS}} 
\nc{\RHS}{\on{RHS}}
\nc{\Pol}{\on{Pol}_{\mu,\hbar}^{\lambda}} 
\nc{\Pols}{\breve{\on{Pol}}_{\mu,\hbar}^{\lambda}} 
\nc{\x}{\varkappa}
\nc{\Rees}{\on{Rees}}

\nc{\consub}{\fext{G_1}{t\mi}}
\nc{\gloop}{\gext{G}{t\mi}} 
\nc{\oo}{\mathcal{O}}
\nc{\tdelta}{{}^{\tau}\Delta} 
\nc{\quo}{\mathcal{X}}
\nc{\ka}{\mathfrak{k}}
\nc{\ga}{\mathfrak{G}} 

\nc{\sdet}{\on{sdet}}
\nc{\sco}{\widehat{S}(u)}
\nc{\hs}{\widehat{s}} 
\nc{\Gr}{\mathsf{Gr}} 

\nc{\Kgr}{\HS_0^\lambda}
\nc{\Kgrr}{\HS_0^{\bar{\lambda}}} 
\nc{\HS}{\mathcal{X}}

\newlength{\mylength}
\DeclareRobustCommand{\SkipTocEntry}[5]{}

\AtBeginDocument{%
   \def\MR#1{}
}

\title[GKLO representations of twisted Yangians]{GKLO representations of twisted Yangians in type $\mathsf{AI}$ and quantizations of symmetric quotients of the affine grassmannian} 

\author[R. Bartlett]{Robin Bartlett}
\address{School of Mathematics Sciences, Queen Mary University of London, London, E1 4NS, United Kingdom} 
\email{\href{mailto:robin.bartlett.math@gmail.com}{robin.bartlett.math@gmail.com}}

\author[T. Prze\'{z}dziecki]{Tomasz Prze\'{z}dziecki}
\address{School of Mathematics, University of Edinburgh, Peter Guthrie Tait Rd, Edinburgh, EH9 3FD, United Kingdom, OrciD: 0000-0001-9700-1007}
\email{\href{mailto:tprzezdz@exseed.ed.ac.uk}{tprzezdz@exseed.ed.ac.uk}}

\author[L. Tappeiner]{Lukas Tappeiner} 
\address{Department of Mathematical Sciences, University of Bath, Claverton Down, Bath BA2 7AY, United Kingdom.}
\email{\href{mailto:lt862@bath.ac.uk}{lt862@bath.ac.uk}} 



\thanks{The first author was partially supported by EPSRC grant EP-R034826-1. The second author was supported by the EPSRC grant No.\ EP/W022834/1 \emph{Kac--Moody quantum symmetric pairs, KLR algebras and generalized Schur--Weyl duality}.}

\begin{document}

\begin{abstract} 
We construct an analogue of Gerasimov--Kharchev--Lebedev--Oblezin \linebreak (GKLO) representations for twisted Yangians of type $\mathsf{AI}$, using the recently found current presentation of these algebras due to Lu, Wang and Zhang. These new representations allow us to define interesting truncations of twisted Yangians, which, in the spirit of Ciccoli--Drinfeld--Gavarini quantum duality, reflect the Poisson geometry of homogeneous spaces. As our main result, we prove that a truncated twisted Yangian quantizes a scheme supported on quotients of transverse slices in the affine grassmannian. 
\end{abstract}

\maketitle

\setcounter{tocdepth}{1}
\tableofcontents


\section{Introduction}

In \cite{GKLO}, Gerasimov, Kharchev, Lebedev and Oblezin (GKLO) constructed a family of infinite-dimensional representations of the Yangian using explicit difference operators. These GKLO representations were later generalized to the shifted Yangian setting and used to produce quantisations of transverse slices to Schubert varieties in the affine Grassmannian \cite{KWWY,BFN}. In this paper we develop an analogue of that picture for twisted Yangians of type $\mathsf{AI}$. Concretely, we  construct explicit GKLO-style representations of these twisted Yangians and realise the resulting truncations as quantisations of loci inside quotients of the thick affine Grassmannian.

Twisted Yangians, introduced by Olshanski \cite{Olsha}, arise in mathematical physics from the Yang--Baxter and reflection equations as algebras controlling the symmetries of integrable systems with boundaries. From an algebraic point of view, they appear naturally in the context of Gelfand--Tsetlin theory for simple Lie algebras of classical types \cite{Molev-GT}. Recently, they have also been realized as degenerations of affine quantum symmetric pair coideal subalgebras, or $\imath$quantum groups \cite{LWZ2}. Shifted versions of twisted Yangians, which we denote by $\dya$ for dominant $\mu \in X_*(T)$, were introduced in \cite{TT} building on the Drinfeld presentation from \cite{LWZ1}.

\subsection*{Setup}
Throughout we take $G = \operatorname{SL}_n$ with diagonal torus $T$ and write $X_*(T) = \operatorname{Hom}(\mathbb{G}_m,T)$. Let $N(T) \subset G$ denote the normaliser of $T$ so that $W = N(T)/T$. We consider the thick affine grassmannian $\operatorname{Gr}^{\operatorname{thick}} =G((z^{-1})) /G[z]$  with the Poisson structure induced from the standard Manin triple for the loop algebra. If $\mathcal{G}_0 \subset G[[z^{-1}]]$ denotes the first congruence subgroup then, for any dominant (always considered relative to the standard upper triangular Borel) $\mu \in X_*(T)$, the $\mathcal{G}_0$-orbits $\operatorname{Gr}_\mu$ through $z^{\mu}$ are known to be Poisson subschemes. We then consider the involution 
$$
\tau: G((z^{-1})) \rightarrow G((z^{-1})), \qquad \tau(g) = g(-z)^t
$$
where $x^t$ denotes the transpose of $x$. Set $K = \lbrace g \in G((z^{-1})) \mid \tau(g) = g^{-1} \rbrace$
and put $K_0 = K \cap \mathcal{G}_0$.

\subsection*{Poisson geometry}
The first part of this paper analyses the Poisson geometry of the quotient spaces $K_0 \backslash \operatorname{Gr}_\mu$. This is done in Section~\ref{affgrassPoisson} and the following summarises our main results.

\begin{thm}\label{thmA}
    Let $K \subset G((z^{-1}))$ denote the subgroup defined by $g(-z)^t = g(z)^{-1}$ and set $K_0 = K \cap \mathcal{G}_0$. Then
    \begin{enumerate}
        \item $K_0$ is a coisotropic subgroup of $\mathcal{G}_0$ and so, for each $\mu \in X_*(T)$, the quotient $K_0 \backslash \operatorname{Gr}_\mu$ (which is representable by an affine scheme) inherits a Poisson structure from  $\operatorname{Gr}_\mu$.
        \item The symplectic leaves inside $K_0 \backslash \operatorname{Gr}_\mu$ are connected components of the subschemes 
        \begin{equation}\label{eq-leaves intro}
         K_0 \backslash \left( K x z^{\eta} \cap \operatorname{Gr}_\mu \right)
        \end{equation}
        where $\eta \in X_*(T)$ and $x \in G$ is such that $x^tx \in N(T)$ represents an involution $w \in W$ with $\eta +w(\eta) \geq 2\mu$ and dominant. The subscheme \eqref{eq-leaves intro} is uniquely determined $\lambda = \eta + w(\eta)$, except when $w$ has no fixed points, in which case there are two such subschemes corresponding to two possible choices of $x$. If $\mu =0$ then \eqref{eq-leaves intro} is connected and hence a single symplectic leaf. 
        \item Write $\mathcal{S}_{2\mu}^\lambda$ for the subscheme \eqref{eq-leaves intro} when $w$ has at least one fixed point and the union of the two distinct such subschemes otherwise. If $\mu =0$ then the ideal of the reduced closed subscheme
        $$
        \mathcal{S}_{2\mu}^{\leq \lambda} = \bigcup_{\gamma\leq \lambda} \mathcal{S}_{2\mu}^{\lambda}
        $$
        inside $K_0 \backslash \operatorname{Gr}_\mu$
        can be described as the radical of an ideal Poisson generated by a set of explicit rational functions defined via trailing principal minors. 
    \end{enumerate}
\end{thm}

A central idea underlying the proof of these results is the existence of an isomorphism
$$
K_0 \backslash \operatorname{Gr}_\mu \cong \operatorname{Gr}_{2\mu}^{\tau=1}
$$
induced via $K_0 x \mapsto \tau(x)x$ for $\tau(g(z)) = g(-z)^t$. This identifies the quotient Poisson structure on the left hand side with that on the right hand side obtained via Dirac reduction (up to a factor of $2$). It also identifies the subschemes $\mathcal{S}_{2\mu}^\lambda$ from Theorem~\ref{thmA} with the $\tau$-fixed points inside the loci $\operatorname{Gr}_{2\mu}^\lambda = \operatorname{Gr}_{2\mu} \cap \operatorname{Gr}^\lambda$ with $\operatorname{Gr}^\lambda \subset \operatorname{Gr}^{\operatorname{thick}}$ the Schubert cell given as the $G[z]$-orbit through $z^\lambda$.

We expect that the restrictions to $\mu=0$ in Theorem~\ref{thmA} are unnecessary. However, we are currently unable to prove this. The most significant obstruction arises in the proof of part (3) of Theorem~\ref{thmA}.

\subsection*{GKLO-representations and quantisation}
 Our second main shows the Poisson structures from Theorem~\ref{thmA} are quantised by   twisted Yangians and their truncations. Write $\dya$ for the $\mathbb{C}[\hbar]$-form of the twisted Yangian shifted by a dominant $\mu \in X_*(T)$.

\begin{thm}\label{thmB}
For each dominant $\mu \in X_*(T)$ there is an isomorphism of Poisson $\mathbb{C}$-algebras
\begin{equation}\label{eq-intro quantise}
\dya / \hbar \dya \cong \mathcal{O}(K_0 \backslash \operatorname{Gr}_\mu). 
\end{equation}
Furthermore,
\begin{enumerate}
    \item For each dominant $\lambda \in X_*(T)$ we define quotients $\dyatr$ of $\dya$ via twisted GKLO-representations. See Theorem~\ref{thm: GKLO}.
    \item If $\mu =0$ then the identification \eqref{eq-intro quantise} induces a Poisson surjection
    $$
    \dyatr / \hbar \dyatr \rightarrow \mathcal{O}(\mathcal{S}_{2\mu}^{\leq -w_0(2\lambda)})
    $$
    which is an isomorphism up to nilpotent elements.
\end{enumerate}
\end{thm}

The twisted GKLO representations from part (1) of Theorem~\ref{thmB} are given in terms of the current, or `new Drinfeld', presentation of the twisted Yangian, found recently by Lu, Wang and Zhang \cite{LWZ1}. The generators commonly denoted as $b_i(u)$ act as sums of localized difference operators, via formulae which somewhat resemble the Gelfand--Tsetlin formulae for quantum symmetric pairs \cite{GavKl, LP}. The shape of our formulae is motivated by the realization of the twisted Yangian in terms of Sklyanin minors. Namely, the corresponding Cartan generators, given by the principal Sklyanin minors, are required to act as \emph{even} polynomials. This restriction determines the correct coefficients on the difference operators corresponding to the $b_i(u)$ generators. The presence of this extra symmetry is a new feature, absent from the original GKLO representations. 

We also expect that the map in part (2) of Theorem~\ref{thmB} is actually an isomorphism, so that $\dya$ directly quantises $\mathcal{S}_{2\mu}^{\leq -w_0(2\lambda)}$, rather than a non-reduced scheme supported on this locus. This could be proved if one knew the Poisson ideal discussed in part (3) of Theorem~\ref{thmA} was reduced. In Conjecture~\ref{conj-reduced} we formulate a conjecture in this direction and show how its validity implies an explicit description of the ideal defining the truncation $\dyatr$ inside $\dya$.

\begin{thm}\label{thmC}
    Assume $\mu  =0$ and that Conjecture~\ref{conj-reduced} holds. Then the surjection in part (2) of Theorem~\ref{thmB} is an isomorphism and $\dyatr$ is the quotient of $\dya$ by a two sided ideal generated by elements $A_i^{(r)}$ for $r \geq r_i := \langle \omega_i, 2\lambda \rangle$ and $B_i^{(r_i+1)}$ for each $1\leq i \leq n-1$.
\end{thm}

\subsection*{A concrete example of Ciccoli--Drinfeld--Gavarini duality}
Theorem~\ref{thmB} can be viewed through the lens of the Drinfeld--Gavarini quantum duality principle \cite{DrinfeldQG, gavarini, gavarini2}. For a Poisson--Lie group $H$ with Lie algebra $\mathfrak{h}$, quantum duality asserts that the semiclassical limit of $U_\hbar(\mathfrak{h})$ is isomorphic, as a Poisson--Hopf algebra, to $\mathcal{O}(H^*)$, the algebra of functions on the dual Poisson--Lie group. Applied to $H=G[z]$ this identifies the semiclassical Yangian with functions on the dual Poisson group $\mathcal{G}_0$ and is the starting point for the results in \cite{KWWY}. Ciccoli--Gavarini \cite{CG06,CG14} extended quantum duality to Poisson homogeneous spaces (quotients by coisotropic subgroups). If $N\subset H$ is coisotropic then the algebra of $N^\perp$-invariant functions on the dual Poisson group corresponds to the semiclassical limit of a (one-sided) coideal subalgebra. In our setting the subgroup $N=G[z]^\tau$ and its orthogonal complement give rise (via this principle) to twisted Yangian-type coideal algebras; Theorem~\ref{thmB} provides a concrete realisation of this correspondence in the split type $\mathsf{AI}$ case.

\begin{rem}
Just before the completion of our paper, another preprint on the same topic was independently released by Lu, Wang and Weekes \cite{LWW-yang}. There is considerable intersection between the two papers. For example, we both introduce GKLO-style representations and study the Poisson geometry of fixed point loci in affine Grassmannian slices. Many aspects of \cite{LWW-yang} are more general, e.g., they also construct GKLO representations for twisted Yangians of type $\mathsf{AIII}$, consider generalized affine Grassmannian slices, and introduce the notion of an iCoulomb branch. 

On the other hand, elements of our approach differ significantly from loc.\ cit.\ and appear to offer more control towards several outstanding open problems. For example, Theorem~\ref{thmB} provides the partial step towards Conjecture 8.13 of loc.\ cit.\ as described in Remark 8.14. In a similar spirit, our Conjecture~\ref{conj-reduced} provides a direct approach towards Conjecture 3.10 of loc.\ cit.\ and presents the possibility of emulating the strategy devised in the untwisted setting in \cite{KMWO}. 
\end{rem}

\addtocontents{toc}{\SkipTocEntry} 
	
\section*{Acknowledgements} 
We would like to thank Gwyn Bellamy, Sasha Shapiro, Lewis Topley and Ryo Yamagishi for helpful discussions and suggestions.

\section{Notation}\label{sec-setup}

\subsection{Basic setup}

	Let $G = \operatorname{SL}_n$ viewed as an algebraic group over $\operatorname{Spec}\mathbb{C}$, with diagonal torus $T$ and upper triangular Borel $B$. Also write $U^{\pm}$ respectively for the upper and lower triangular unipotent subgroups in $G$. Set $W = N(T)/T$  for $N(T) \subset G$ the normaliser of $T$ and write $w_0 \in W$ for the longest element.
    
  Set $X_*(T) = \operatorname{Hom}(\mathbb{G}_m,T)$ and $X^*(T) = \operatorname{Hom}(T,\mathbb{G}_m)$, which are dual via the evaluation pairing 
	$$
	\langle-, - \rangle: X_*(T) \times X^*(T) \rightarrow \mathbb{Z}, \qquad (\beta^\vee, \lambda) \mapsto \langle \beta^\vee, \lambda \rangle = \beta^\vee(\lambda).
	$$
We use analogous notations for the diagonal torus $\widetilde{T} \subset \operatorname{GL}_n$, with canonical maps 
$X_*(T) \hookrightarrow X_*(\widetilde{T})$ and $X^*(\widetilde{T}) \twoheadrightarrow X^*(T)$. 

	We will use the index sets $ \indx = \{1, \cdots, n-1\}$ and $\tilde{\indx} = \{1, \cdots, n\}$. Let $(a_{ij})_{i,j\in \indx}$ be the Cartan matrix associated with $G$. Given $i \in \tilde{\indx}$, we write $\epsilon_i^\vee \in X^*(\widetilde{T})$ for the character 
	$$
	\operatorname{diag}(t_1,\ldots,t_n)  \mapsto t_i
	$$
	and $\epsilon_i$ for its dual under $\langle -,-\rangle$. 
Let $\Delta^+$ be the set of all positive roots of $G$ (relative to $B$), with simple roots $\alpha^\vee_i = \epsilon_i^\vee - \epsilon_{i+1}^\vee$ $(i \in \indx)$. We write $\omega_i \in X_*(T)\otimes_\Z \Q$ for the fundamental coweights, dual to $\alpha_i^\vee$, and $\omega_i^\vee \in X^*(T) $ for the fundamental weights. Let $X_*(T)^+ \subset X_*(T)$ denote the set of dominant coweights relative to $B$, so $\lambda \in X_*(T)^+$ iff $\langle \alpha_i^\vee, \lambda \rangle \geq 0$ for each $i \in \indx$ and write $\mu \leq\lambda$ if $\lambda - \mu$ is a sum of positive coroots. . 

\subsection{Affine grassmannians}\label{sec-setupgrass} Write $G((z^{-1}))$, $G[z]$, and $G[z,z^{-1}]$ for the ind-group schemes over $\operatorname{Spec}\mathbb{C}$ with $A$-valued points given respectively by $G(A((z^{-1})), G(A[z])$ and $G(A[z,z^{-1}])$. We use similar notation when $G$ is replaced by another affine group scheme over $\operatorname{Spec}\mathbb{C}$. We then consider the fpqc quotients
$$
\operatorname{Gr}^{\operatorname{thick}} = G((z^{-1}))/ G[z], \qquad \operatorname{Gr}^{\operatorname{thin}} = G[z,z^{-1}]/G[z]
$$
the first of which is representable by a scheme, and the latter by an ind-scheme. For any coweight $\lambda \in X_*(T)$ we write $z^\lambda \in \operatorname{Gr}^{\operatorname{thick}}$ for the image of $\lambda(z) \in G((z^{-1}))$ and if $\lambda$ is dominant we write $\operatorname{Gr}^{\lambda}$ for the $G[z]$-orbit through this point, with closure $\operatorname{Gr}^{\leq \lambda}$. We also consider the subgroup 
$$
\mathcal{G}_0 \subset G((z^{-1}))
$$
whose $A$-points, for any $\mathbb{C}$-algebra $A$, consist of matrices in $1 + z^{-1}\operatorname{Mat}(A[[z^-1]])$. Recall that if $\mathcal{U}^{\pm}_0 = \mathcal{G}_0 \cap U^{\pm}((z^{-1}))$ and $\mathcal{T}_0 = \mathcal{G}_0 \cap T((z^{-1}))$ then multiplication defines an isomorphism
\begin{equation}\label{eq-Gauss factorisation}
\mathcal{U}^+_0 \times \mathcal{T}_0 \times \mathcal{U}^-_0 \rightarrow \mathcal{G}_0
\end{equation}
(indeed multiplication $U^+ \times T \times U^- \rightarrow G$ is known to be an open immersion whose image is the open locus defined by the non-vanishing of the principal minors). For $\mu \in X_*(T)^+$ set $\operatorname{Gr}_\mu \subset \operatorname{Gr}^{\operatorname{thick}}$ equal the $\mathcal{G}_0$-orbit through $z^\mu$. Finally, for any $1 \leq i \leq n$ and any pair of $i$-tuples $I,J \subset \lbrace 1,\ldots,n\rbrace$ we write 
$$
\Delta_{IJ} \in \mathcal{O}(\mathcal{G}_0)[[z^{-1}]]
$$
for the series valued function whose value on $g \in \mathcal{G}_0$ is its $IJ$-th minor. Write $\Delta_{IJ}^{(r)} \in \mathcal{O}(\mathcal{G}_0)$ for the coefficient of $z^{-r}$ in this series.
\subsection{Commutators}\label{sec-setupcomm}
	We  use the following notation for commutators: $[a,b] = ab - ba$ and $[a,b]_+ = ab + ba$.

\section{Symmetric quotients of the affine grassmannian}\label{affgrassPoisson}

\subsection{Poisson structures on loop groups}

Following \cite{KWWY} we equip $G((z^{-1}))$ with the Poisson structure induced by the Manin triple $(\mathfrak{g}((z^{-1})),\mathfrak{g}[z],z^{-1}\mathfrak{g}[[z^{-1}]])$, with the pairing  
$$(x,y) = \operatorname{Res}_{z=0} \operatorname{Trace}(xy)
$$
i.e., the coefficient of $z^{-1}$ in $\operatorname{Trace}(xy)$. In \cite[Proposition 2.13]{KWWY} (specialised to the case $G = \operatorname{SL}_n$) the resulting Poisson bracket $\lbrace-,-\rbrace$ is computed as:

\begin{lem}\label{lem-kwwypoissonstruc}
	Recall the series valued functions $\Delta_{IJ} = \Delta_{IJ}(z)$ from Section~\ref{sec-setupgrass}. Then
	$$
	\lbrace \Delta_{IJ}(u),\Delta_{KL}(v) \rbrace = \frac{1}{u-v}  \sum_{1 \leq p,q \leq n} \left( \epsilon^{J,L}_{pq} \Delta_{IJ^{(p \leadsto q)}}(u) \Delta_{KL^{(q \leadsto p)}}(v)-  \epsilon^{I,K}_{qp} \Delta_{I^{(q \leadsto p)}J}(u) \Delta_{K^{(p \leadsto q)}L}(v)  \right)
	$$
	where 
	\begin{itemize}
			\item $J^{(p \leadsto q)}$ denotes the tuple with $p \in J$ replaced by $q$. Likewise for $L^{(q \leadsto p)}, I^{(q \leadsto p)},$ and $K^{(p \leadsto q)}$.
		\item $\epsilon^{J,L}_{pq} =0$ if $p \not\in J$ or $q \not\in L$ and otherwise equals $\pm 1$ according to the sign of the permutations reordering $J^{(p \leadsto q)}$ and $L^{(q \leadsto p)}$ into ascending order.
	\end{itemize}
\end{lem} 

Both $G[z]$ and $\mathcal{G}_0$ appear as Poisson subgroups of $G((z^{-1}))$ and the quotient $\operatorname{Gr}^{\operatorname{thick}} = G((z^{-1})) / G[z]$ inherits a Poisson structure from that on $G((z^{-1}))$.

\begin{lem}\label{lem-lemma about U's}
	Recall $\mathcal{U}_0^- = \mathcal{G}_0 \cap U^-((z^{-1}))$. For $\mu \in X_*(T)$ dominant, set $\mathcal{U}^{-,\mu} = \mathcal{G}_0 \cap z^\mu G[z] z^{-\mu}$ and $\mathcal{U}^{-}_\mu = z^\mu \mathcal{U}_0^- z^{-\mu}$. Then both are subgroups of $\mathcal{U}_0^-$ and:
	\begin{enumerate}
		\item $\mathcal{U}^{-,\mu}$ is a coisotropic subgroup in $\mathcal{G}_0$, and so $\mathcal{G}_0 \rightarrow \mathcal{G}_0 / \mathcal{U}^{-,\mu}$ is a Poisson quotient.
		\item The quotient map $\mathcal{U}^-_\mu \rightarrow \mathcal{U}^-_0 / \mathcal{U}^{-,\mu}$ is an isomorphism.
	\end{enumerate}
	Identical assertions hold with $\mu$ anti-dominant and each $\mathcal{U}^-_0$ replaced by $\mathcal{U}^+_0 = \mathcal{G}_0 \cap U^+((z^{-1}))$.
\end{lem}
\begin{proof}
	For part (1) one shows that the orthogonal complement of $\operatorname{Lie}\mathcal{U}^{-,\mu}$ in $\mathfrak{g}[z]$ is a Lie subalgebra, which is an easy computation. Part (2) follows since multiplication $\mathcal{U}^{+}_{\mu} \times \mathcal{U}^{+,\mu} \rightarrow \mathcal{U}_0^+$ is an isomorphism.
\end{proof}
Recall $\operatorname{Gr}_\mu \subset \operatorname{Gr}^{\operatorname{thick}}$ is the $\mathcal{G}_0$-orbit through $z^\mu$ for dominant $\mu \in X_*(T)$. Lemma~\ref{lem-lemma about U's} furnishes two descriptions of this locally closed subscheme. Firstly, the orbit map induces an isomorphism $\mathcal{G}_0 / \mathcal{U}^{-,\mu} \xrightarrow{\sim} \operatorname{Gr}_\mu$ which, by part (1) of Lemma~\ref{lem-lemma about U's}, endows $\operatorname{Gr}_\mu$ with a natural Poisson structure. On the other hand, part (2) of Lemma~\ref{lem-lemma about U's} combined with the factorisation \eqref{eq-Gauss factorisation} shows how acting on the base point in $\operatorname{Gr}^{\operatorname{thick}}$ gives an isomorphism $\mathcal{W}_\mu := \mathcal{U}^+_0 \mathcal{T}_0 u^\mu \mathcal{U}^- \xrightarrow{\sim} \operatorname{Gr}_\mu$. 

In what follows it will also be useful to consider a variant of this construction giving an isomorphism
\begin{equation}\label{eq-isom}
\mathcal{U}^{+,-\mu} \backslash \mathcal{G}_0 / \mathcal{U}^{-,\mu} \rightarrow \mathcal{U}^+_{-\mu} \times \mathcal{T}_0 \times \mathcal{U}^-_\mu \xrightarrow{ g \mapsto z^\mu g z^\mu}\mathcal{W}_{2\mu}
\end{equation}
where the first map comes from part (2) of Lemma~\ref{lem-lemma about U's}. We claim \eqref{eq-isom} is Poisson for the quotient Poisson structure on the source and that on the target induced by the isomorphism $\mathcal{W}_{2\mu} \cong \operatorname{Gr}_{2\mu}$. To see this note the composite $\mathcal{G}_0 \rightarrow\mathcal{U}^{+,-\mu} \backslash \mathcal{G}_0 / \mathcal{U}^{-,\mu} \rightarrow \mathcal{W}_{2\mu}$ is exactly the shift morphism denoted $\iota_{2\mu,-\mu,-\mu}$ in \cite[\S5.9]{Toda}. That this is Poisson follows from \cite[Theorem 5.15]{Toda}.
\subsection{Symmetric quotients and fixed points}
Consider the anti-involution $\tau(g) = g^{t}(-z)$ on $G((z^{-1}))$ and set
$$
K = \lbrace g \in G((z^{-1})) \mid \tau(g) = g^{-1} \rbrace
$$ 
This can be interpreted as the loop group associated to the special unitary group over $\mathbb{C}((z^{-1}))$ respecting the Hermitian  form $(x,y) = x^t(z) y(-z)$ on $\mathbb{C}((z^{-1}))^n$. Thus, \cite[Theorem 0.1]{PR08} ensures $K$ is connected. The group $K$ acts on $\operatorname{Gr}^{\operatorname{thick}}$ via left multiplication.

\begin{rem}
Previous work \cite{Na04,CY03} consider the action on $\operatorname{Gr}^{\operatorname{thick}}$ of fixed points in $G((z^{-1}))$ under involutions like $g \mapsto g^{-t}$ on $\operatorname{Gr}$. While there are a number of similarities, the geometry of this action is different to ours in several significant  ways.
\end{rem}

We avoid discussion of the quotient stack $K\backslash \operatorname{Gr}^{\operatorname{thick}}$ and instead restrict attention to $K_0 \backslash \operatorname{Gr}_\mu$ for $K_0 = K \cap \mathcal{G}_0$.  The latter are represented by schemes, since $g \mapsto \tau(g) g$ induces a monomorphism
\begin{equation}\label{eq-Psi}
	\Psi: K_0 \backslash \operatorname{Gr}_\mu = K_0 \backslash \mathcal{G}_0 / \mathcal{U}^{-,\mu} \rightarrow \mathcal{U}^{+,-\mu} \backslash \mathcal{G}_0 /\mathcal{U}^{-,\mu} \xrightarrow{\eqref{eq-isom}} \mathcal{W}_{2\mu}^{\tau =1}
\end{equation}
which is easily checked, e.g. by comparing Hilbert series of the tangent spaces at the identity, graded via loop rotations, to be surjective and hence an isomorphism. 

\begin{lem}
	The subgroup $K_0 \subset \mathcal{G}_0$ is coisotropic. 
\end{lem}
\begin{proof}
	As in Lemma~\ref{lem-lemma about U's}, this follows since the orthogonal complement of $\operatorname{Lie}K_0 $
	 in $\mathfrak{g}[z]$ is $\operatorname{Lie}(K \cap G[z])$ which is a Lie subaglebra.
\end{proof}

Thus $K_0 \backslash \operatorname{Gr}_\mu = K_0 \backslash \mathcal{G}_0 / \mathcal{U}^{-,\mu}$ has a natural Poisson structure and \eqref{eq-Psi} transfers this structure to $\mathcal{W}_{2\mu}^{\tau =1 } $. On the other hand, $\left(  \mathcal{U}^{+,-\mu} \backslash \mathcal{G}_0 /\mathcal{U}^{-,\mu} \right)^{\tau =1} \cong \mathcal{W}_{2\mu}^{\tau = 1}$ has an intrinsic Poisson structure induced from that on $\operatorname{Gr}_{2\mu} \cong \mathcal{W}_{2\mu}$ via Dirac reduction (see, for example, \cite[\S 5.4.3]{LGPV}, \cite[\S 2]{To23}). As explained in \cite[Proposition 5.36]{LGPV}, the corresponding Poisson bracket $\lbrace-,-\rbrace_\tau$ is given by 
\begin{equation}\label{eq-Dirac bracket}
\lbrace F,G \rbrace_\tau = \frac{1}{2} \left( \lbrace \widetilde{F}, \widetilde{G} \rbrace + \lbrace \tau^* \widetilde{F}, \widetilde{G} \rbrace \right)	
\end{equation}
where $\lbrace-,-\rbrace$ denotes the Poisson bracket on $\mathcal{W}_{2\mu}$ and $\widetilde{F},\widetilde{G} \in \mathcal{O}(\mathcal{W}_{2\mu})$ are lifts of $F,G \in \mathcal{O}(\mathcal{W}_{2\mu}^{\tau =1})$. In fact, these two Poisson structures on $\mathcal{W}_{2\mu}^{\tau = 1}$ coincide up to a multiple of $2$, as follows from \cite[Theorem 5.9]{Xu}. More precisely, one applies loc.\ cit.\ (which, while written in a finite dimensional setting, goes through immediately in our loop setup) in the case $\mu =0$, and then deduces the claim for $\mu >0$ using functoriality of Dirac reduction. 

\subsection{Symplectic leaves and their closures}

\begin{defi}
	For dominant $\lambda \geq 2\mu$ write 
	$$
	\mathcal{S}_{2\mu}^{ \lambda} := \Psi^{-1}(\mathcal{W}_{2\mu}^{\tau =1} \cap G[z]z^\lambda G[z] )
	$$
	for $\Psi$ as in \eqref{eq-Psi} and consider the closed subset $\mathcal{S}_{2\mu}^{\leq \lambda} = \bigcup_{\gamma \leq \lambda} \mathcal{S}_{2\mu}^\gamma$ equipped with the reduced scheme structure.
\end{defi}

Recall from \cite[Theorem 2.5]{KWWY} that the intersections $\mathcal{W}_{2\mu} \cap G[z] z^\lambda G[z]$ are symplectic leaves inside $\mathcal{W}_{2\mu}$. These exhaust all symplectic leaves meeting $G[z,z^{-1}]$. It follows from \cite[Proposition 5.26]{LGPV}  that the connected components of $\mathcal{W}_{2\mu}^{\tau =1} \cap G[z] z^\lambda G[z]$ are then the symplectic leaves in $\mathcal{W}_{2\mu}^{\tau =1}$. The same is therefore true of the $\mathcal{S}_{2\mu}^\lambda$. 

\begin{pro}\label{prop-eta+w(eta)}
	If $\mathcal{S}_{2\mu}^\lambda \neq \emptyset$ then $\lambda = \eta + w(\eta)$ for some $\eta \in X_*(T)^+$ and some involution $w \in W$ (i.e. $w^2=1$). If one can take $w$ with at least one fixed point then
	$$
	\mathcal{S}_{2\mu}^\lambda = K_0 \backslash \left( K x u^\eta \cap \operatorname{Gr}_\mu \right)
	$$
	for any $x \in G$ with $\tau(x)x \in N(T)$ lifting $w$. Otherwise there are $x^\pm \in G$ with $\tau(x^\pm)x^\pm \in N(T)$ lifting $w$ so that $K_0 \backslash \left( K x^\pm u^\eta \cap \operatorname{Gr}_\mu \right)$ are distinct---in this case $\mathcal{S}_{2\mu}^\lambda$ is the disjoint union of these two intersections.
\end{pro}

Note that if $\lambda = \sum_{1 \leq i \leq n} \lambda_i \epsilon_i$ then $\lambda = \eta + w(\eta)$ as in the Proposition~\ref{prop-eta+w(eta)} if and only if the number of $\lambda_i$ equal any given odd number is even. The second case of Proposition~\ref{prop-eta+w(eta)} occurs when every $\lambda_i$ is odd (and is therefore only possible when $n$ is even and $\geq 4$).

\begin{proof}
		Consideration of the map \eqref{eq-Psi} shows that $K_0 q \in K_0 \backslash \operatorname{Gr}_\mu$ lies inside $\mathcal{S}_{2\mu}^\lambda$ if and only if $\tau(q_0)q_0 \in G[z] z^\lambda G[z]$ for any $q_0 \in G((z^{-1}))$ representing $q \in \operatorname{Gr}_\mu$. Notice this implies $K_0 q \cap \operatorname{Gr}^{\operatorname{thin}} \neq \emptyset$.
		
		On the other hand, \cite[Theorem 5.2]{DF03} shows that any $q \in \operatorname{Gr}^{\operatorname{thin}}$ can be represented by $q_0 \in G[z,z^{-1}]$ with $\tau(q_0)q_0 = \dot{w}z^\lambda$ where
		\begin{itemize}
			\item $\lambda = (\lambda_1,\ldots, \lambda_n) \in X_*(T)^+$,
			\item $\dot{w} \in N(T) \cap G^\tau$ represents an involution $w \in W$ with $w(\lambda) = \lambda$ and whose fixed points are exactly the $i \in \lbrace 1,\ldots,n \rbrace$  with $\lambda_i \in 2\mathbb{Z}$ and with $w(\lambda) = \lambda$.
		\end{itemize}
		This immediately implies $\lambda = \eta + w(\eta)$ for some $\eta \in X_*(T)^+$. Now every element in $ N(T) \cap G^\tau$ can be expressed as $\tau(x) x$ for some $x \in G$. Applying this to $(-1)^{\eta_\lambda} \dot{w}$ allows us to write $\dot{w} z^\lambda = \tau(xz^{\eta}) x z^{\eta}$. We conclude that the $K$-orbit through $q$ equals $K x z^\eta$. Consequently,
		$$
		\mathcal{S}_{2\mu}^\lambda = K_0 \backslash \left( \bigcup Kx u^{\eta} \cap \operatorname{Gr}_\mu \right)
		$$
		with the union running over $x \in G$ with $\tau(x)x \in N(T)$ lifting $w$. It only remains to determine when two such $Kxu^\eta$ coincide, and this reduces to a description of when the $K \cap G$-orbits through the class of $x$ inside $G/ P_\eta$ coincide for $P_\eta \subset G$ the parabolic subgroup stabilising $u^\eta \in \operatorname{Gr}^{\operatorname{thick}}$. Such a description is given in \cite[Lemma 10.3.1]{RS90} and the proposition follows.
\end{proof}

\begin{lem}
    If $\mu =0$ then each subscheme $K_0 \backslash (K x z^\eta \cap \operatorname{Gr}_\mu)$ as in Proposition~\ref{prop-eta+w(eta)} is connected, and hence a symplectic leaf in $K_0 \backslash \operatorname{Gr}_\mu$.
\end{lem}
\begin{proof}
    It suffices to show $K x z^\eta \cap \operatorname{Gr}_\mu$ is connected. If it is empty then we are done so assume not. As $\operatorname{Gr}_0$ is open in $\operatorname{Gr}^{\operatorname{thick}}$  this intersection is open inside the $K$-orbit through $x z^\eta$. But $K$ is a connected group and so this orbit is irreducible. The same must then be true of any non-empty open subset. 
\end{proof}
For later use we record the following simple combinatorial consequence of the constraint on $\lambda$ given in Proposition~\ref{prop-eta+w(eta)}.

\begin{cor}\label{cor-constraint on the r_i's}
    Suppose $\lambda = \sum_{1 \leq i \leq n} \lambda_i \epsilon_i$ can be expressed as $\eta + w(\eta)$ for an involution $w \in W$. If
    $$
    r_i := \langle \omega_i, -w_0(\lambda) \rangle = -(\lambda_{n-i+1}+\ldots+\lambda_n)
    $$
    then $r_i$ odd implies $\lambda_{n-i} = \lambda_{n-i+1}$ are both odd, and so $r_{i+1}$ is even.
\end{cor}

In general, we do not know whether $\mathcal{S}_{2\mu}^{\leq \lambda}$ coincides with the closure of $\mathcal{S}_{2\mu}^\lambda$ in $K_0 \backslash \operatorname{Gr}_\mu$ whenever $\mathcal{S}_{2\mu}^\lambda$ is non-empty. However, we are able to show this holds in the simplest case, which is when $\lambda$ is even, i.e. lies in $2X_*(T)^+$.

\begin{thm}\label{thm-even leaves}
	If $\lambda \in 2X_*(T)^+$ then $\mathcal{S}_{2\mu}^{ \leq \lambda}$ equals the closure of $\mathcal{S}_{2\mu}^\lambda$ in $K_0 \backslash \operatorname{Gr}_\mu$ and coincides with the scheme theoretic image of $\operatorname{Gr}_{\mu} \cap \operatorname{Gr}^{\leq \eta} \rightarrow K_0 \backslash \operatorname{Gr}_\mu$ (recall the notation from Section~\ref{sec-setupgrass}).
\end{thm}

\begin{proof}
	Write $\overline{\mathcal{S}}_{2\mu}^\lambda$ for the closure of $\mathcal{S}_{2\mu}^\lambda$ in $K_0 \backslash \operatorname{Gr}_\mu$ and set $\lambda = 2\eta$. We then consider the $K^0$-orbit $\mathcal{O}^0$ through $u^\eta$ where $K^0 := K \cap G[z]$. We first claim that the closure of $\mathcal{O}^0$ in $\operatorname{Gr}^{\operatorname{thick}}$ equals $\operatorname{Gr}^{\leq \eta}$. Certainly this closure is contained in $\operatorname{Gr}^{\leq \eta}$ and, since the latter is irreducible, equality follows if $\mathcal{O}^0$ and $\operatorname{Gr}^{\leq \eta}$ have the same dimension, i.e. if  $\operatorname{dim} \mathcal{O}^0 =   \sum_{\alpha^\vee >0} \langle \alpha^\vee, \eta \rangle$. Since the orbit map identifies  $\mathcal{O}^0 \cong K^0 / (K^0 \cap z^\eta G[z] z^{-\eta})$ this dimension can be computed as 
	$$
	\operatorname{dim} \operatorname{Lie}K^0 / \left( \operatorname{Lie}K^0 \cap z^\eta \mathfrak{g}[z] z^{-\eta} \right)
	$$
	But $\operatorname{Lie}K^0 = \lbrace x \in \mathfrak{g}[z] \mid \tau(x) = -x \rbrace$ is spanned by $z \operatorname{Lie}T[z^2]$ together with the elements $z^i (x_{\alpha^\vee} - (-1)^r x_{-\alpha^\vee})$ for all $i \geq 0$ and $\alpha^\vee >0$. As $\eta$ is dominant it follows that the above quotient is spanned by the images of $z^i (x_{\alpha^\vee} - (-1)^r x_{-\alpha^\vee})$ for $0 \leq i < \langle \alpha^\vee, \eta \rangle$ and $\alpha^\vee >0$. This proves the claim.
	
	Next we show that $\operatorname{Gr}^{\leq \lambda} \cap \operatorname{Gr}_\mu$ equals the closure in $\operatorname{Gr}_\mu$ of $\mathcal{O}^0 \cap \operatorname{Gr}_\mu$. Since $\mathcal{O}^0$ is dense in $\operatorname{Gr}^{\leq \eta}$ by the previous paragraph and $z^\mu \operatorname{Gr}_0$ is open in $\operatorname{Gr}^{\operatorname{thick}}$ with $z^\mu \operatorname{Gr}_0 \cap \operatorname{Gr}^{\leq \eta}$ non-empty, we deduce that $\mathcal{O}^0\cap z^\mu\operatorname{Gr}_0$ is non-empty and has closure $\operatorname{Gr}^{\leq \eta} \cap \operatorname{Gr}_{\mu}$ inside $\operatorname{Gr}_\mu$. To replace $z^\mu \operatorname{Gr}_0$ with $\operatorname{Gr}_\mu$ in this assertion notice that
	$$
	z^\mu \operatorname{Gr}_0 \cong z^\mu \mathcal{U}^+_0 z^{-\mu}  \mathcal{T}_0 z^\mu \mathcal{U}^-_0 \cong  \left( z^\mu \mathcal{U}^+_0 z^{-\mu} \cap G[z] \right) \times \mathcal{W}_\mu
	$$
	where the second isomorphism uses Lemma~\ref{lem-lemma about U's} to identify $\mathcal{U}^+_0 = \mathcal{U}^{+,-\mu} \times  \mathcal{U}^{+}_{-\mu}$, and so identify $z^\mu\mathcal{U}^+_0 z^{-\mu} = \left( z^\mu \mathcal{U}^+_0 z^{-\mu} \cap G[z]\right) \times \mathcal{U}^+_0$. Let $p: z^\mu \operatorname{Gr}_0 \rightarrow \operatorname{Gr}_\mu$ be the resulting projection. Since  $\operatorname{Gr}^{\leq \eta}$  is $G[z]$-stable we have $p^{-1}(\operatorname{Gr}^{\leq \eta} \cap \operatorname{Gr}_\mu) = \operatorname{Gr}^{\leq \eta} \cap z^\mu \operatorname{Gr}_0$. Now suppose $Z \subset \operatorname{Gr}^{\leq \eta} \cap \operatorname{Gr}_\mu$  is closed and contains $\mathcal{O}^0 \cap \operatorname{Gr}_\mu$. By the first assertion of this paragraph it follows that $p^{-1}(Z) =  \operatorname{Gr}^{\leq \eta} \cap z^\mu \operatorname{Gr}_0$ and so $Z = \overline{\mathcal{O}}^0 \cap \operatorname{Gr}_\mu$. This shows $\operatorname{Gr}^{\leq \eta} \cap \operatorname{Gr}_\mu$ is the closure of $\mathcal{O}^0 \cap \operatorname{Gr}_\mu$ as required.
	
	Proposition~\ref{prop-eta+w(eta)} shows $\mathcal{S}_{2\mu}^\lambda = K_0 \backslash \left( K z^\eta \cap \operatorname{Gr}_\mu \right)$ when $\lambda = 2\eta$. The well-known fact that multiplication $\mathcal{G}_0 \times G[z] \rightarrow G((z^{-1}))$ is an open immersion implies the same for $K_0 \times K^0 \rightarrow K$. As a consequence, the image of $\mathcal{O}^0 \cap \operatorname{Gr}_\mu \rightarrow K_0 \backslash \operatorname{Gr}_\mu$ is dense inside $\mathcal{S}_{2\mu}^\lambda$. The previous paragraph therefore shows that $\operatorname{Gr}^{\leq \lambda} \cap \operatorname{Gr}_\mu \rightarrow K_0 \backslash \operatorname{Gr}_\mu$ factors through $\overline{\mathcal{S}}_{2\mu}^\lambda$ with dense image. It follows that $\overline{\mathcal{S}}_{2\mu}^\lambda$ is the scheme theoretic image of this morphism.
	
	Now suppose $\gamma \leq \lambda$ with $\mathcal{S}_{2\mu}^\gamma$ non-empty. Combining Proposition~\ref{prop-eta+w(eta)} with the openness of multiplication $K_0 \times K^0 \rightarrow K$ shows that, as $\nu$ runs over coweights with $\gamma = \nu + w(\nu)$ with $w \in W$ an involution represented by $\tau(g)g$ for $g \in G$, the images of $K^0 x u^{\nu} \cap \operatorname{Gr}_\mu$ in $K_0 \backslash \operatorname{Gr}_\mu$ are each dense in a union of connected components of $\mathcal{S}_{2\mu}^\gamma$. Since $\gamma \leq \lambda$, we have $\nu \leq \eta$. Thus, any such $K^0 x u^{\nu} \cap \operatorname{Gr}_\mu$ is contained inside $\operatorname{Gr}^{\leq \eta} \cap \operatorname{Gr}_\mu$. We conclude that $\mathcal{S}_{2\mu}^\gamma \subset \overline{\mathcal{S}}_{2\mu}^\lambda$. This finishes the proof since clearly $\overline{\mathcal{S}}_{2\mu}^\lambda \subset \bigcup_{\gamma\leq \lambda} \mathcal{S}_{2\mu}^{\gamma}$.	\end{proof}

    For the remaining $\lambda \geq 2\mu$, i.e. those not in $2X_*(T)^+$ but of the form $\eta + w(\eta)$ for an involution $w$, we expect that $\mathcal{S}_{2\mu}^\lambda$ is non-empty. Indeed, computations made when $\mu =0$ suggest that $\mathcal{S}_{2\mu}^\lambda$ is irreducible (except in those cases described in Proposition~\ref{prop-eta+w(eta)}) of dimension
    $$
    \sum_{1 \leq i \leq n-1} 2 \lfloor{ \frac{n_{i}}{2} \rfloor} 
    $$
    for $n_{i}$ defined by $\lambda = 2\mu + \sum_{1 \leq i \leq n-1} n_{i} \alpha_i$.
    
    \subsection{Ideal generators}\label{sec-ideal generators}
    
    Our goal in this section is to describe Poisson generators of the ideal of $\mathcal{S}_{2\mu}^{ \leq\lambda}$. We are currently only able to do this when $\mu =0$ and we impose this restriction throughout this section. In particular, $\Psi: K_0\backslash \mathcal{G}_0 \rightarrow \mathcal{G}_0^{\tau=1}$ from \eqref{eq-Psi} is simply the map $K_0x_0 \mapsto \tau(x_0)x_0$.

    \begin{notat}\label{notation-tau minors}
        For $i$-tuples $I,J \subset \lbrace 1,\ldots,n\rbrace$ recall $\Delta_{IJ} \in \mathcal{O}(\mathcal{G}_0)$ is defined in Section~\ref{sec-setupgrass}. Set
        $$
        \Delta_{IJ}^{\tau} :=  \Delta_{IJ} \circ \Psi \in \mathcal{O}(K_0 \backslash \mathcal{G}_0)
        $$
        and $\Delta_{IJ}^{\tau,(r)}$ for the coefficient of $z^{-r}$ in this series. Thus $\Delta_{IJ}^\tau(K_0x_0) = \Delta_{IJ}(\tau(x)x)$.
    \end{notat}
    
    Of particular relevance will be the functions
    \begin{equation}
    A_i^{(r)}:=  \Delta_{\lbrace n-i+1,\ldots,n \rbrace,\lbrace n-i+1,\ldots,n \rbrace}^{\tau,(r)}, \qquad B_i^{(r)}:= \Delta_{\lbrace n-i+1,\ldots,n \rbrace,\lbrace n-i,n-i+2,\ldots,n \rbrace}^{\tau,(r)}
	\end{equation}
    lying inside $\mathcal{O}(K_0\backslash \mathcal{G}_0)[[z^{-1}]]$.

    \begin{pro}\label{prop-nonvanishing}
        Suppose $\mathcal{S} \subset \mathcal{S}_{0}^\lambda$ is a symplectic leaf with $\lambda = \sum_{1 \leq i \leq n} \lambda_i \epsilon_i \in X_*(T)^+$ not necessarily in $2X_*(T)^+$ and set 
        $$
        r_i =  \langle \omega_i, -w_0 \lambda \rangle = -(\lambda_n + \ldots + \lambda_{n-i+1})
        $$
        Then, as functions on $\mathcal{S}$, $A_{i+1}^{(r_{i+1})} \neq 0$ if $r_{i+1}$ is even, and $B_{i+1}^{(r_{i+1})} \neq 0$  if $r_{i+1}$ is odd.
    \end{pro}

Before giving the proof recall that if $g \in G[z]z^\lambda G[z]$ then, for each pair of $i$-tuples $I,J \subset \lbrace 1,\ldots,n\rbrace$, the minor
$\Delta_{IJ}(g)$
has $z$-adic valuation $\geq \langle \omega_i, -w_0(\lambda)\rangle$. Furthermore, for each $1 \leq i \leq n$ there is an $I,J$ for which this is an equality. In the proof of Proposition~\ref{prop-nonvanishing} we will need the following more specific version of this assertion:

    \begin{lem}\label{lem-next minor is off by only two}
	Let $\lambda = \sum_{1\leq i \leq n} \lambda_i \epsilon_i$ and  suppose $g \in G[z]z^\lambda G[z]$ is such that \linebreak $\Delta_{\lbrace n-i+1,\ldots,n \rbrace, \lbrace n-i+1,\ldots,n \rbrace}(g)$ has $z$-adic valuation $\lambda_n + \ldots+ \lambda_{n-i+1}$. Then there are $1 \leq r,s \leq n-i$ for which
	$$
	\Delta_{\lbrace r, n-i+1,\ldots,n \rbrace, \lbrace s, n-i+1,\ldots,n \rbrace}(g)
	$$
	has $z$-adic valuation $\lambda_n + \ldots + \lambda_{n-i}$.
\end{lem}
\begin{proof}
	The hypothesis on $g$ says that if $g=\left( \begin{smallmatrix}
		X & Y \\ Z & W  
	\end{smallmatrix} \right)$ with $W$ an $i-1$ by $i-1$ matrix then $\operatorname{det}(W)$ has $z$-adic valuation $\lambda_n+ \ldots+ \lambda_{n-i+1}$. Since the maximal minors of $Z$ and $Y$ have $z$-adic valuation at least that of $\operatorname{det}W$ we can find matrices $Y^*, Z^*$ so that 
	$$
	g' := \begin{pmatrix}
		1 & Y^* \\ 0 & 1 
	\end{pmatrix} \begin{pmatrix}
	X & Y \\ Z & W 
	\end{pmatrix} \begin{pmatrix}
	1 & 0 \\ Z^* & 1 
	\end{pmatrix} = \begin{pmatrix}
	X' & 0 \\ 0 & W
	\end{pmatrix}
	$$
	Note $\Delta_{\lbrace r, n-i+1,\ldots,n \rbrace, \lbrace s, n-i+1,\ldots,n \rbrace}(g') = \Delta_{\lbrace r, n-i+1,\ldots,n \rbrace, \lbrace s, n-i+1,\ldots,n \rbrace}(g)$ for any $1 \leq r,s \leq n-i$, as can be seen by considering the action of $ \left( \begin{smallmatrix}
		1 & Z^* \\ 0 & 1	\end{smallmatrix} \right)$ and $ \left( \begin{smallmatrix}
	1 & 0 \\ Y^* & 1
	\end{smallmatrix} \right) $ on the relevant element inside $\bigwedge^{i} \mathbb{C}^n$. On the other hand, the fact that $ \left( \begin{smallmatrix}
	X' & 0 \\ 0 & W	\end{smallmatrix} \right) \in G[z] z^\lambda G[z]$ ensures $\Delta_{\lbrace r, n-i+1,\ldots,n \rbrace, \lbrace s, n-i+1,\ldots,n \rbrace}(g')$ has $z$-adic valuation $\lambda_n + \ldots + \lambda_{n-i}$ for some $1 \leq r,s \leq n-i$. 
\end{proof}

    \begin{proof}[Proof of Proposition~\ref{prop-nonvanishing}]

        Argue by induction on $i+1$. If $r_{i}$ is even (or $i =0$) then we can assume the existence of $K_0 x_0 \in K_0 \backslash \mathcal{G}_0$ with $A_{i}^{(r_{i})}(K_0x_0) \neq 0$. Applying Lemma~\ref{lem-next minor is off by only two} therefore produces $1 \leq r,s \leq n-i$ with
        $$
        \Delta^{\tau,(r_{i+1})}_{\lbrace r, n-i+1,\ldots,n \rbrace, \lbrace s, n-i+1,\ldots,n \rbrace} 
        $$
        non-vanishing on $K_0x_0$. To prove the proposition in this case we will vary the point $K_0x_0$ inside $\mathcal{S}$ by flowing along suitable Hamiltonian vector fields. 
        
        More precisely, we will use the following observation: If $f \in \mathcal{O} (K_0 \backslash \mathcal{G}_0)$ and $x \in \mathcal{S}_{0}^{\lambda}$ is a closed point with
	$$
	\lbrace f,g \rbrace_\tau(x) \neq 0
	$$
	for some $g \in \mathcal{O} (K_0 \backslash \mathcal{G}_0)$ then there exists a closed point $y \in \mathcal{S}_{0}^{\lambda}$ with $f(y) \neq 0$. Indeed, $\mathcal{S}_{0}^\lambda$ being a symplectic leaf means that the tangent space $T_x \mathcal{S}_{0}^\lambda$ is spanned by the value at $x$ of the Hamiltonian vector fields on $K_0 \backslash \mathcal{G}_0$, while $\lbrace f, g \rbrace_\tau(x)$ is by definition the value of $f$ on the Hamiltonian vector field associated to $g$ at $x$. Since $\mathcal{S}_{0}^\lambda$ is smooth, and hence reduced, it follows that the vanishing locus of $f$ has codimension $1$ in $ \mathcal{S}_{0}^\lambda$.

    \subsubsection*{Step 1} First, we claim $x_0$ can be chosen so that at least one of $r$ and $s$ equals $1$. Suppose for a contradiction that this is not the case. Then \eqref{eq-Dirac bracket} and  Lemma~\ref{lem-kwwypoissonstruc} together give
        $$
        \begin{aligned}
        \lbrace \Delta^{\tau,(r_{i+1})}_{\lbrace n-i, n-i+1,\ldots,n \rbrace, \lbrace s, n-i+1,\ldots,n \rbrace},& \Delta_{n-i,r}^{\tau,(1)} \rbrace_\tau = \\&= \frac{1}{2}\lbrace \Delta^{(r_{i+1})}_{\lbrace n-i, n-i+1,\ldots,n \rbrace, \lbrace s, n-i+1,\ldots,n \rbrace}, \Delta_{n-i,r}^{(1)} \rbrace \circ \Psi  \\
        &\quad -\frac{1}{2} \lbrace \Delta^{(r_{i+1})}_{\lbrace n-i, n-i+1,\ldots,n \rbrace, \lbrace s, n-i+1,\ldots,n \rbrace}, \Delta_{r,n-i}^{(1)} \rbrace \circ \Psi \\&= \frac{1}{2}\Delta^{\tau,(r_{i+1})}_{\lbrace r , n-i+1,\ldots,n \rbrace, \lbrace s, n-i+1,\ldots,n \rbrace}
        \end{aligned}
        $$
        as functions on $\mathcal{S}$ (we remind the reader here that $\lbrace-,-\rbrace_\tau$ denotes the Poisson bracket on $K_0\backslash \mathcal{G}_0$ pulled back via $\Psi$ from that on $\mathcal{G}_0^{\tau=1}$ obtained via Dirac reduction from the bracket $\lbrace-,-\rbrace$ on $\mathcal{G}_0$). But by assumption $\Delta^{\tau,(r_{i+1})}_{\lbrace r , n-i+1,\ldots,n \rbrace, \lbrace s, n-i+1,\ldots,n \rbrace}$ does not vanish on $\mathcal{S}$ and so the argument of the previous paragraph produces a point in $\mathcal{S}$ on which $\Delta^{\tau,(r_{i+1})}_{\lbrace n-i, n-i+1,\ldots,n \rbrace, \lbrace s, n-i+1,\ldots,n \rbrace}$ does not vanish.

        \subsubsection*{Step 2} If $r_{i+1}$ is even then 
        $$
        \tau^* \Delta^{(r_{i+1})}_{\lbrace n-i,\ldots,n \rbrace, \lbrace n-i, \ldots,n \rbrace} = \Delta^{(r_{i+1})}_{\lbrace n-i,\ldots,n \rbrace, \lbrace n-i, n-i+1,\ldots,n \rbrace}$$
        inside $\mathcal{O}(\mathcal{G}_0)$. This, together with  \eqref{eq-Dirac bracket} and Lemma~\ref{lem-kwwypoissonstruc} combine to give 
        $$
        \begin{aligned}
        \lbrace \Delta^{\tau,(r_{i+1})}_{\lbrace n-i, \ldots,n \rbrace, \lbrace n-i, \ldots,n \rbrace}, \Delta_{s,n-i}^{\tau,(1)} \rbrace_\tau &= \lbrace \Delta^{(r_{i+1})}_{\lbrace n-i, \ldots,n \rbrace, \lbrace n-i, n-i+1,\ldots,n \rbrace}, \Delta_{s,n-i}^{(1)} \rbrace \circ \Psi \\ &= \Delta^{(\tau,r_{i+1})}_{\lbrace n-i , \ldots,n \rbrace, \lbrace s, n-i+1,\ldots,n \rbrace}
        \end{aligned}
        $$
        as functions on $\mathcal{S}$. By Step 1 we know this function is non-vanishing, and so we produce a point of $\mathcal{S}$ on which $\Delta^{\tau,(r_{i+1})}_{\lbrace n-i, n-i+1,\ldots,n \rbrace, \lbrace n-i, n-i+1,\ldots,n \rbrace}$ does not vanish.

        \subsubsection*{Step 3} If instead $r_{i+1}$ is odd then we can assume $s \neq n-i-1$ (since otherwise we are done) and $s \neq n-i$ (since then $\Delta^{\tau,(r_{i+1})}_{\lbrace n-i, n-i+1,\ldots,n \rbrace, \lbrace s, n-i+1,\ldots,n \rbrace} =0$). Using \eqref{eq-Dirac bracket} together with Lemma~\ref{lem-kwwypoissonstruc} then gives the identity
        $$
        \begin{aligned}
        \lbrace & \Delta^{\tau.(r_{i+1})}_{\lbrace n-i, \ldots,n \rbrace, \lbrace n-i-1,n-i+1,\ldots,n \rbrace}, \Delta_{s,n-i-1}^{(1)} \rbrace_\tau  \\
        &= \lbrace \Delta^{(r_{i+1})}_{\lbrace n-i, \ldots,n \rbrace, \lbrace n-i-1,n-i+1,\ldots,n \rbrace}, \Delta_{s,n-i-1}^{(1)} \rbrace\circ \Psi - \lbrace \Delta^{(r_{i+1})}_{\lbrace n-i, \ldots,n \rbrace, \lbrace n-i-1,n-i+1,\ldots,n \rbrace}, \Delta_{n-i-1,s}^{(1)} \rbrace \circ \Psi
        \\
        &= - \Delta^{\tau,(r_{i+1})}_{\lbrace n-i,\ldots,n \rbrace, \lbrace s,n-i+1,\ldots,n \rbrace} 
        \end{aligned}
        $$
        Again, by Step 1 we know this function is non-vanishing, and so we produce a point of $\mathcal{S}$ on which $\Delta^{\tau,(r_{i+1})}_{\lbrace n-i, n-i+1,\ldots,n \rbrace, \lbrace n-i-1, n-i+1,\ldots,n \rbrace}$ does not vanish.

        \subsubsection*{Step 4} We conclude with the case $r_i$ is odd. By assumption $\mathcal{S}_{0}^\lambda$ is non-empty and so Corollary~\ref{cor-constraint on the r_i's} forces $r_{i+1}$ to be even and $\lambda_{n-i+1} = \lambda_{n-i}$ to be odd. The inductive hypothesis gives $K_0x \in K_0 \backslash \mathcal{G}_0$ for which the series valued function
        $$
        \Delta_{\lbrace n-i+1,\ldots,n \rbrace,\lbrace n-i,n-i+2,\ldots,n \rbrace}^{\tau}
        $$
        has $z$-adic valuation $r_i$. On the other hand, the Desnanot–Jacobi formula \cite{VV23} gives
        \begin{equation}\label{eq-alice in wonderland equation}
        \begin{aligned}
            \Delta_{\lbrace n-i,\ldots,n \rbrace,\lbrace n-i,\ldots,n \rbrace}^\tau \Delta_{\lbrace n-i+2,\ldots,n\rbrace, \lbrace n-i+2,\ldots,n \rbrace}^\tau= \alpha -\beta
        \end{aligned}
        \end{equation}
        for
        $$
        \begin{aligned}
        \alpha &= \Delta_{\lbrace n-i+1,\ldots,n \rbrace,\lbrace n-i+1,\ldots,n \rbrace}^\tau \Delta_{\lbrace n-i,n-i+2,\ldots,n\rbrace, \lbrace n-i,n-i+2,\ldots,n \rbrace}^\tau, \\
        \beta &= \Delta_{\lbrace n-i,n-i+2,\ldots,n \rbrace,\lbrace n-i+1,\ldots,n \rbrace}^\tau \Delta_{\lbrace n-i+1,\ldots,n\rbrace, \lbrace n-i,n-i+2,\ldots,n \rbrace}^\tau
        \end{aligned}
        $$
        A priori the two minors in $\alpha$ have $z$-adic valuation $\geq r_{i}$. In fact this inequality is strict because $r_{i}$ is odd and so $\tau$-invariance forces the coefficient of $z^{-r_{i}}$ to vanish. Therefore, the right hand side of \eqref{eq-alice in wonderland equation} has $z$-adic valuation $2r_i$. But recall $\lambda_{n-i+1} =  \lambda_{n-i}$ and so $2r_i = r_{i+1} + r_{i-1}$. We conclude that $\Delta^\tau_{\lbrace n-i,\ldots,n \rbrace,\lbrace n-i,\ldots,n \rbrace}$ and $\Delta^\tau_{\lbrace n-i+2,\ldots,n\rbrace, \lbrace n-i+2,\ldots,n \rbrace}$, which respectively have $z$-adic valuations $\geq r_{i+1}$ and $\geq r_{i-1}$, in fact have exact $z$-adic valuations $r_{i+1}$ and $r_{i-1}$. In particular, $A_{i+1}^{(r_{i+1})} \neq 0$ on $\mathcal{S}$ and we are done. 
    \end{proof}

    \begin{cor}\label{cor-ideal generators}
		Assume $\lambda \in X_*(T)$ can be written as $\eta +w(\eta)$ for an involution $w \in W$. Then the ideal defining $\mathcal{S}_{0}^{ \leq \lambda}$ inside $\mathcal{O}(K_0 \backslash \mathcal{G}_0)$ is the radical of the ideal Poisson generated by the functions
        \begin{itemize}
            \item $A_{i}^{(r)}$
		for $1 \leq i \leq n-1$ and $r > \langle \omega_i, -w_0\lambda \rangle$.
        \end{itemize}
	\end{cor}
    \begin{proof}
        Write $J_{0}^\lambda$ for ideal Poisson generated by the elements described, and $V(J_0^\lambda) \subset K_0 \backslash \operatorname{Gr}_0$ for its vanishing locus. Clearly $\mathcal{S}^{\leq \lambda}_0 \subset V(J_0^\lambda)$ because if $I,J$ are $i$-tuples then $\Delta_{IJ}^{\tau,(r)}$ vanishes on $\mathcal{S}_0^\lambda$ whenever  $r > \langle \omega_i,-w_0\lambda \rangle$.

        For the opposite inclusion we claim that $J_0^\lambda$ contains all $\Delta_{ij}^{\tau,(r)}$ with $1 \leq i,j \leq n$ and sufficiently large $r$. This ensures that the image of  $V(J_0^\lambda)$ under $\Psi$ lies inside $\operatorname{Gr}^{\operatorname{thin}}$, and so that $V(J_0^\lambda)$ is a union of the symplectic leaves $\mathcal{S} \subset \mathcal{S}_{0}^{\gamma}$ for varying $\gamma$. Granting this, the corollary follows from Proposition~\ref{prop-nonvanishing}. Indeed, if $\gamma = \sum_{1\leq i \leq n} \gamma_i \epsilon_i> \lambda =\sum_{1\leq i \leq n} \lambda_i \epsilon_i$ then there is a smallest $i$ with $r_i' := \langle \omega_i, -w_0\delta \rangle >\langle \omega_i, -w_0\lambda \rangle := r_i$. If $r_i'$ is even then Proposition~\ref{prop-nonvanishing} shows $A_i^{(r_i')} \neq 0$ on $\mathcal{S}$ and so $\mathcal{S} \not\subset V(J_0^\lambda)$.  If $r_i'$ is odd then Corollary~\ref{cor-constraint on the r_i's} implies $r_{i+1}' = r_{i}' - \gamma_{n-i} = r_{i}'-\gamma_{n-i+1}$ is even. On the other hand, since $r_{i-1}' = r_{i-1}$ and $r_i'>r_i$ we must have $-\gamma_{n-i+1} > - \lambda_{n-i+1}$. Thus, $r_{i+1}' > r_{i+1}$ and Proposition~\ref{prop-nonvanishing} shows $A_i^{(r_{i+1}')} \neq 0$ on $\mathcal{S}$, hence $\mathcal{S} \not\subset V(J_0^\lambda)$ and we are done.

        It only remains to prove our claim. If $i,j \neq 1$ and $r$ is  even then Lemma~\ref{lem-kwwypoissonstruc} combined with \eqref{eq-Dirac bracket} gives $\lbrace \Delta_{1,1}^{\tau,(r)}, \Delta_{1,j}^{\tau,(1)} \rbrace_\tau = \Delta_{1,j}^{\tau,(r)}$ and
        $$
         \lbrace \Delta_{1,j}^{\tau,(r)}, \Delta_{i,1}^{\tau,(1)} \rbrace_\tau = \frac{1}{2}\lbrace \Delta_{1,j}^{(r)},\Delta_{i,1}^{(1)} \rbrace \circ \Psi+ \frac{1}{2}\lbrace \Delta_{j,1}^{(r)},\Delta_{i,1}^{(1)} \rbrace \circ \Psi= \frac{1}{2}\left(  \delta_{ij} \Delta_{11}^{\tau,(r)} - \Delta_{ij}^{\tau,(r)}  \right)
        $$
        This shows $J_0^\lambda$ contains $\Delta_{ij}^{\tau,(r)}$ for all even $r >r_1$. On the other hand, if $r$ is even then \eqref{eq-Dirac bracket} combined with Lemma~\ref{lem-kwwypoissonstruc} gives
        $$
        \lbrace \Delta_{ij}^{\tau,(r)}, \Delta_{ll}^{\tau,(2)} \rbrace_{\tau} = \lbrace \Delta_{ij}^{(r)}, \Delta_{ll}^{(2)} \rbrace \circ \Psi =  -\Delta_{il}^{\tau,(r+1)} \delta_{lj} +\Delta_{lj}^{\tau,(r+1)} \delta_{il} -\Delta_{il}^{\tau,(r)} \Delta_{lj}^{\tau,(1)} + \Delta_{lj}^{\tau,(r)} \Delta_{il}^{\tau,(1)}
        $$
        We've already seen that the final two terms on the right lie in $J_0^\lambda$. Therefore $\Delta^{\tau,(r+1)}_{il}\delta_{lj} - \Delta_{lj}^{\tau,(r+1)} \delta_{il} \in J_0^\lambda$. It follows that $\Delta_{il}^{\tau,(r+1)} \in J_0^\lambda$ if $i \neq l$ and $r > r_1$ is even. Since $\Delta_{ii}^{\tau,(r+1)} =0$ for $r$ even this finishes the proof.
    \end{proof}

    We emphasise that it is essential to consider the radical of the ideal in Corollary~\ref{cor-ideal generators}. For example, suppose $\lambda \in 2X_*(T)^+$. Then $B_1^{(r_1+1)}$ vanishes on $\mathcal{S}_0^\lambda$ but the ideal Poisson generated by the $A_i^{(r)}$'s for $r> r_i$ only consists of  functions with loop grading $\geq r_1+2$ (since the bracket has degree $-1$ for this grading). We believe that this non-reducedness can be eliminated by adding the $B_i^{(r_i+1)}$'s to the Poisson generating set whenever $r_i$ is even.

    \begin{conjecture}\label{conj-reduced}
        The ideal of $\mathcal{O}(K_0 \backslash \mathcal{G}_0)$ which is Poisson generated by the functions
        \begin{itemize}
            \item $A_{i}^{(r)}$
		for $1 \leq i \leq n-1$ and $r > \langle \omega_i, -w_0\lambda \rangle$,
        \item $B_i^{(r_i+1)}$ whenever $r_i = \langle \omega_i, -w_0(\lambda) \rangle$ is even
        \end{itemize} 
        is reduced, and hence is the ideal defining $\mathcal{S}_0^{ \leq \lambda}$ inside $K_0 \backslash\mathcal{G}_0$. \end{conjecture}

\section{Shifted Yangians}
	
	Below we recall the definitions and main properties of shifted Yangians for $\sln$.

\subsection{Drinfeld presentation} 

Given $\mu \in X_*(T)^+$, let 
\[
H_i(u) = u^{\ai(\mu)} + \h \sum_{r= -\ai(\mu)}^\infty H_i^{(r)} u^{-r-1}. 
\]
We call $\underline{H}_i(u) = \sum_{r=0}^\infty H_{i}^{(r)} u^{-r-1}$ the \emph{principal part} of $H_i(u)$ (and use analogous notation for principal parts of other series). Also let 
\[
E_i(u) = \h \sum_{r=0}^\infty E_i^{(r)} u^{-r-1}, \qquad  F_i(u) = \h \sum_{r=0}^\infty F_i^{(r)} u^{-r-1}. 
\] 

The \emph{$\mu$-shifted Yangian} $Y_{\mu,\h}$ associated to $\sln$ is the $\C[\h]$-algebra generated by $H_i^{(r)}, E_i^{(s)}$, $F_i^{(s)}$ (with $i \in \indx, r \geq -\ai(\mu), s \in \N$), subject to the following defining relations: 
\begin{align*}
[H_i(u), E_j(v)] &= - \half a_{ij} \h \frac{[H_i(u), E_j(u) - E_j(v)]_+}{u-v}, \\ 
[H_i(u), F_j(v)] &= \half a_{ij} \h \frac{[H_i(u), F_j(u) - F_j(v)]_+}{u-v}, \\  
[E_i(u), F_j(v)] &= \delta_{ij} \h \frac{\underline{H}_i(u) - \underline{H}_i(v)}{u-v}, \\ 
[E_i(u), E_i(v)] &= -  \h \frac{(E_i(u) - E_i(v))^2}{u-v}, \\ 
[F_i(u), F_i(v)] &=  \h \frac{(F_i(u) - F_i(v))^2}{u-v}, 
\end{align*}
\begin{align*}
[E_i(u), E_j(v)] &= \half  \h \frac{[E_i(u), E_j(u) - E_j(v)]_+}{u-v} - \frac{[E_i^{(0)}, E_j(u) - E_j(v)]}{u-v} \quad (|j-i| =1), \\ 
[F_i(u), F_j(v)] &= - \half  \h \frac{[F_i(u), F_j(u) - F_j(v)]_+}{u-v} - \frac{[F_i^{(0)}, F_j(u) - F_j(v)]}{u-v} \ \ (|j-i| =1),  \end{align*}
\begin{align*}
\Sym_{u_1,u_2} [E_i(u_1), [E_i(u_2), E_j(v)]] &= 0 \qquad (|j-i| =1),  \\ 
\Sym_{u_1,u_2} [E_i(u_1), [E_i(u_2), E_j(v)]] &= 0 \qquad (|j-i| =1),  \\
[E_i(u), E_j(v)] = [F_i(u), F_j(v)] &= 0 \qquad (a_{ij} = 0).  
\end{align*} 
We denote the specialization to $\hbar = 1$, i.e., $Y_{\mu,\h}/(\h-1)Y_{\mu,\h}$, as $Y_\mu$. If $\mu = 0$, we also abbreviate $Y_{\h} = Y_{0, \h}$ and $Y = Y_0$.  

\subsection{PBW theorem} 

For each positive root $\beta^\vee =  \alpha_{j} + \alpha_{j+1} + \cdots + \alpha_{i}$ $(j \leq i)$, and $r \geq 0$, set 
\begin{align}
E_{\beta^\vee}^{(r)} &= [[ \cdots [E_j^{(r)}, E_{j+1}^{(0)}], \cdots, E_{i-1}^{(0)}], E_i^{(0)} ], \\
F_{\beta^\vee}^{(r)} &= [F_{i}^{(0)}, [F_{i-1}^{(0)}, \cdots, [F_{j+1}^{(0)}, F_{j}^{(r)}] \cdots ]]. 
\end{align}
We call elements $E_{\beta^\vee}^{(r)}, F_{\beta^\vee}^{(r)}$ and $H_i^{(s_i)}$, as $\beta^\vee$ ranges over $\Delta^+$, $i$ over $\indx$, $r \geq 0$ and $s_i \geq -\ai(\mu)$, PBW variables. Fix any total ordering on the set of PBW variables. 

\Thm[{\cite[Theorem 2.55]{FT}}] \label{thm: Toda PBW}
Ordered PBW monomials in the PBW variables form a basis of $Y_{\mu,\h}$ as a free $\C[\h]$-module, as well as a basis of $Y_\mu$ as a free $\C$-module. 
\enthm 

As a consequence of Theorem \ref{thm: Toda PBW} (c.f. \cite[Corollary 3.16]{Toda}), given anti-dominant coweights $\mu_1,\mu_2$, the shift homomorphisms 
\eq \label{eq: shift homo untwisted} 
\iota(\mu,\mu_1,\mu_2): Y_{\mu} \rightarrow Y_{\mu+\mu_1+\mu_2} 
\eneq
defined by
\[
H_i^{(r)} \mapsto H_i^{(r - \ai(\mu_1+\mu_2))},  \qquad 
E_i^{(r)} \mapsto E_i^{(r - \ai(\mu_1))},  \qquad 
F_i^{(r)} \mapsto F_i^{(r - \ai(\mu_2))}. 
\] 
are injective. In particular, if $\mu_1 + \mu_2 = - \mu$, then $\iota(\mu,\mu_1,\mu_2)$ yields an embedding $Y_{\mu} \hookrightarrow Y$. 


\subsection{Canonical filtration} 

For the purpose of quantizing slices in the affine Grassmannian, one requires a certain subalgebra of $Y_{\mu,\h}$. 
Namely, let $\Yu$ be the $\C[\h]$-subalgebra of $Y_{\mu,\h}$ generated by 
\[
\{\h E_{\beta^\vee}^{(r)}\}_{\beta^\vee \in \Delta^+}^{r \geq 0} \cup 
\{\h F_{\beta^\vee}^{(r)}\}_{\beta^\vee \in \Delta^+}^{r \geq 0} \cup 
\{\h H_{i}^{(s_i)}\}_{i \in \indx}^{s_i \geq - \ai(\mu)}. 
\]
This subalgebra also admits an alternative description as a Rees algebra, which we now recall. 
By \cite[\S 5.4]{Toda}, given a pair of coweights $\mu_1, \mu_2$ such that $\mu_1 + \mu_2 = \mu$, there is a filtration $F^\bullet_{\mu_1,\mu_2}$ on $Y_\mu$ determined by 
\[
\deg E_{\beta^\vee}^{(r)} = \beta^\vee(\mu_1) + r +1, \qquad \deg F_{\beta^\vee}^{(r)} = \beta^\vee(\mu_2) + r+1, \qquad \deg H_{i}^{(r)} = \alpha_i^\vee(\mu) + r+1. 
\] 
According to \cite{Toda}, this defines an algebra filtration, is independent of the choice of PBW variables and their ordering, and the corresponding Rees algebra $\Rees^{F^\bullet_{\mu_1,\mu_2}} Y_\mu$ is independent of the choice of coweights $\mu_1,\mu_2$. 

\Thm[{\cite[Theorem 2.57]{FT}}]
For any $\mu \in X_*(T)^+$, there is a canonical $\C[\h]$-algebra isomorphism
\[
\Yu \cong \Rees^{F^\bullet_{\mu_1,\mu_2}} Y_\mu. 
\]
\enthm 

Passing to the semi-classical limit, there are canonical isomorphisms
\[
\Yu / \h \Yu \cong \Rees^{F^\bullet_{\mu_1,\mu_2}} Y_\mu / \h \Rees^{F^\bullet_{\mu_1,\mu_2}} Y_\mu \cong \gr^{F^\bullet_{\mu_1,\mu_2}} Y_\mu. 
\]
Moreover, by \cite[Proposition 5.7]{Toda}, $\gr^{F^\bullet_{\mu_1,\mu_2}} Y_\mu$ is a commutative algebra. 
It is also naturally a Poisson algebra, with the Poisson bracket given by 
\[
\{ a, b \} = \hbar\mi [\hat{a}, \hat{b}] \quad \mbox{mod} \ \h\Yu, 
\]
for any lifts $\hat{a}, \hat{b}$. 
Finally, we remark that the shift homomorphisms \eqref{eq: shift homo untwisted} induce monomorphisms of Poisson algebras
\[
\Yu/\h \Yu \hookrightarrow \mathbf{Y}_{\mu + \mu_1 + \mu_2} / \h \mathbf{Y}_{\mu + \mu_1 + \mu_2}. 
\]

\subsection{Quantum duality} 
\label{sec: quantum duality}

When $\mu = 0$, the algebra $\mathbf{Y} = \mathbf{Y}_0$ can be identified with the Drinfeld--Gavarini dual of $Y_{\h}$. 
Let us first recall the general definition of this notion. 

Let $\mathfrak{a}$ be a Lie algebra over $\C$, and suppose that $A$ is a deformation--quantization of the Hopf algebra $U(\mathfrak{a})$, i.e., $A$ is a Hopf algebra over $\C[\h]$, and there is an isomorphism of Hopf algebras $A/\h A \cong U(\mathfrak{a})$. Let $\Delta$ and $\epsilon$ be the coproduct and counit of $A$, respectively. For $m \geq 0$, define inductively $\Delta^m \colon A \to A^{\otimes n}$ by 
\[
\Delta^0 = \epsilon, \qquad \Delta^1 = \on{id}, \qquad \Delta^m = (\Delta \otimes \on{id}^{\otimes (m-2)}) \circ \Delta^{m-1}. 
\] 
Also define 
\[
\delta_m \colon A \to A^{\otimes m}, \qquad \delta_m = (\on{id} - \epsilon)^{\otimes m} \circ \Delta^m. 
\]
The \emph{Drinfeld--Gavarini dual} $A'$ of $A$ is the following sub-Hopf algebra: 
\[
A' = \{ a \in A \mid \delta_m(a) \in \h^m A^{\otimes m} \mbox{ for all } m \in \N \}. 
\] 
The importance of $A'$ derives from the fact that, according to the quantum duality principle \cite[Theorem 1.6]{gavarini}, $A'$ is a deformation--quantization of the coordinate ring of an algebraic group $G_{\mathfrak{a}^*}$ associated to the dual Lie algebra $\mathfrak{a}^*$. 

Let us now return to the specific setting of the Yangian $Y_\h$. It is a graded Hopf algebra, with $\deg(\h) = 1$ and $\deg(x^{(r)}) = r$ for $x = H_i, E_{\beta^\vee}, F_{\beta^\vee}$. The description of the Hopf algebra structure can be found in, e.g., \cite[(A.16)--(A.17)]{FT}. The Yangian is a deformation--quantization of the universal enveloping algebra of the current Lie algebra, i.e., $Y_\h / \h Y_\h \cong U(\sln[t])$. The latter is an isomorphism of graded Hopf algebras if $U(\sln[t])$ is endowed with the loop grading.
\Thm[{\cite[Corollary A.22]{FT}}]
There is a canonical $\C[\h]$-algebra isomorphism 
\[
Y_\h' \cong \mathbf{Y}. 
\]
\enthm

\subsection{RTT presentation} 

The \emph{RTT Yangian} $\Yartt$ for $\mathfrak{gl}_n$ is the algebra generated by $t_{ij}^{(k)}$, where $k \geq 1$ and $1 \leq i,j \leq n$, subject to the relations 
\[
(u-v) [t_{ij}(u),t_{kl}(v)] =\hbar (t_{kj}(u)t_{il}(v) - t_{kj}(v)t_{il}(u)), 
\]
where 
\[
t_{ij}(u) = \delta_{ij} + \hbar\sum_{r >0} t_{ij}^{(r)} u^{-r}. 
\]
Writing 
\[
T(u) = \sum_{i,j=1}^n e_{ij} \otimes t_{ij}(u) , \qquad P = \sum_{i,j=1}^n e_{ij} \otimes e_{ji}    \qquad R(u) = 1 - Pu\mi,
\]
where $e_{ij}$ are the usual matrix units, 
these relations can be rewritten as
\[
R(u-v)T_1(u)T_2(v) = T_2(v)T_1(u)R(u-v). 
\]

For any formal series $f(u) \in 1 + \frac{\h}{u} \C[\h] [[u^{-1}]]$, the assignment 
\eq \label{eq: rtt auto} 
T(u) \mapsto f(u) T(u)
\eneq
defines an algebra automorphism of $\Yartt$. The RTT Yangian $\sYartt$ for $\sln$ is the $\C[\h]$-subalgebra of $\Yartt$ consisting of all the elements fixed under all automorphisms \eqref{eq: rtt auto}. 

Let $\bYartt$ be the $\C[\h]$-subalgebra of $\Yartt$ generated by $\{ \h t_{ij}^{(r)}\}_{1 \leq i,j \leq n}^{r \geq 1}$, and let $\bsYartt = \bYartt \cap \sYartt$. 

\Thm \label{thm: Yrtt=Y}
There is a $\C[\h]$-algebra isomorphism $Y_{\h} \cong \sYartt$, which restricts to an isomorphism $\mathbf{Y} \cong \bsYartt$. 
\enthm 

\Proof
The isomorphism is constructed by Gauss decomposition, see, e.g., \cite[\S 1.11]{Molev-book}. For the latter statement, see \cite[Proposition 2.34]{FT}. 
\enproof

In light of Theorem \ref{thm: Yrtt=Y}, we will, from now on, just use the notations $Y_\h$ and $\mathbf{Y}$. 


\section{Shifted twisted Yangians}\label{sec-shifted twisted Yangians}

Below we recall the definition of shifted twisted Yangians of type $\mathsf{AI}$ (in the current presentation) and their main properties. 
While we are principally interested in the twisted Yangians associated to $\sln$, we will also need their $\gln$-version. The latter has the advantage that it is easier to verify that the $\gln$-version of the relations are preserved by the twisted GKLO homomorphism defined in \S \ref{sec: tGKLO}.

\subsection{Generators and relations} 

The Drinfeld presentations of the twisted Yangians of split type associated to $\gln$ and $\sln$ (i.e., type $\mathsf{AI}$)  were first found in \cite[Theorem 5.1, Theorem 5.3]{LWZ1} via Gauss decomposition. A generalization to all simple split types, excluding $\mathsf{G}_2$, was established in \cite[\S 4]{LWZ2}. Presentations in terms of generating currents, which are more convenient for our purposes, can be found in \cite[Remark 4.12]{LWZ1} and \cite[\S 4.3]{LWZ2}. 
Below, in the $\gln$ case, we use a presentation with a slightly different choice of Cartan generators than in \cite[Theorem 5.1]{LWZ1}, which can be found in \cite[Proposition 5.3]{LPTTW}.  
The definition of shifted twisted Yangians, based on the Drinfeld presentation, was first given in \cite[\S 3.2]{TT}. Later, an equivalent parabolic presentation was given in \cite[\S 8]{LPTTW}. Below, we pursue the Drinfeld approach. 

Given $\mu \in X_*(T)^+$, 
let 
\[
h_i(u) = u^{-2\ei(\mu)} + \hbar \sum_{r=2\ei(\mu)}^\infty h_{i}^{(r)} u^{-r-1}, \qquad 
b_j(u) = \hbar \sum_{r=0}^\infty b_{j}^{(r)} u^{-r-1}, 
\] 
for $i \in \eindx$ and $j \in \indx$. 
We also set 
\[
\tih_i(u) = (h_i(u))^{-1}, \qquad z_i(u) = \tih_i(u - \half \hbar) h_{i+1}(u). 
\]

\Defi 
\label{def: tYangian} 
The $\mu$-shifted twisted Yangian $\Tya$ of split type associated to $\mathfrak{gl}_n$ is the $\hcor$-algebra generated by $h_{i}^{(r)}, \tih_{i}^{(s)}, b_{j}^{(t)}$ (with $i \in \eindx$, $j \in \indx$, $r \geq 2\ei(\mu)$, $s \geq -2\ei(\mu) $ and $t \in \N$), subject to the following defining relations: 
\begin{align}
    z_i(u) =& \, z_i(-u), \label{rel: h1} \\
    [h_{i}(u),h_{j}(v)]=&\, 0, \label{rel: hh} \\
    [h_{i}(u),b_{j}(v)]=&\,\frac{\delta_{ij}\hbar}{u-v+\tfrac{1}{2}\hbar}\big(b_{i}(u+\tfrac{1}{2}\hbar)-b_{i}(v)\big)h_{i}(u) \label{rel: hb}\\
    &\, +\frac{\delta_{ij}\hbar}{u+v+\tfrac{1}{2}\hbar} h_{i}(u)\big(b_{i}(v)-b_{i}(-u-\tfrac{1}{2}\hbar)\big) \\
    &\, + \frac{\delta_{i,j+1}\hbar}{u+v}h_{i}(u)\big(b_{j}(-u)-b_{j}(v)\big)\label{h2b n2}\\ 
    &\, +\frac{\delta_{i,j+1}\hbar}{u-v} \big(b_{j}(v)-b_{j}(u)\big)h_{i}(u), \\
    [b_i(u),b_{j}(v)] =&\, 0 \quad \textnormal{if } |i-j|>1, \label{rel: bibj} \\
    [b_{i}(u),b_{i}(v)]=& \, \frac{\hbar}{v-u}\big(b_{i}(v)-b_{i}(u)\big)^2 +\frac{\hbar}{u+v}\big(\uz_{i}(v)-\uz_{i}(u)\big), \label{eq: bbz} \\ \label{rel: bibi1}
    (u-v)[b_{i}(u), b_{i+1}(v)]
=&\, -\half \hbar \big(b_i(u)b_{i+1}(v)+b_{i+1}(v)b_i(u)\big) \\
&\, + [b_{i}^{(0)}, b_{i+1}(v)]+[b_{i+1}^{(0)},b_i(u)], 
\end{align}
and
\begin{align} \label{rel: Serre} 
&(u_1+u_2)\Sym_{u_1,u_2}\Big[b_{i}(u_1),\big[b_{i}(u_2), b_{j}(t) \big]\Big] =
\\
&= 4\h \sum_{u = u_1,u_2}  \frac{u(t-\h )\, \uz_i(u) b_j(t)- u(t+\h )\, b_j(t)\uz_i(u) }{4u^2-\h^2}
\end{align}
if $|i-j|=1$. 
\enDefi 

We denote the specialization to $\hbar = 1$ as $\Tyar$. The usual unshifted twisted Yangian is $\Ty = {}^{\mathsf{tw}}\widetilde{{Y}}_{0,\hbar}$. The $\mu$-shifted twisted Yangian $\sTya$ of split type associated to $\mathfrak{sl}_n$ is the subalgebra of $\Tya$ generated by $b_{j}^{(s)}$ and $z_{i}^{(r)}$ ($i \in \tilde{\indx}$,\ $s \geq 0$,\ $r \geq -2\ai(\mu)$). 
We also write $\sTyar$ for the specialization of $\sTya$ at $\hbar = 1$, and $\sTy = {}^{\mathsf{tw}}{{Y}}_{0,\hbar}$.

Consider briefly the special case when $\mu = 0$. Let $Z(\Ty)$ denote the centre of $\Ty$. 

\Pro
There is an isomorphism 
\[
\Ty \cong \sTy \otimes_{\C[\h]} Z(\Ty). 
\]
Moreover, $Z(\Ty) = \C[\h][c_1, c_3, \cdots]$, 
where
\eq \label{eq: central series}
c(u) = 1 + \hbar\sum_{r \geq 0} c_r u^{-r-1} := h_1(u)h_2(u-\hbar) \cdots h_n(u-\hbar(n-1)). 
\eneq
\enpro

\Proof
See \cite[Theorems 2.8.2, 2.9.2]{Molev-book}.  
\enproof

\subsection{PBW theorem} 

For each positive root $\beta^\vee =  \alpha_{j} + \alpha_{j+1} + \cdots + \alpha_{i}$ $(j \leq i)$, and $r \geq 0$, set
\[ 
b_{\beta^\vee}^{(r)} = [b_{i}^{(0)}, [b_{i-1}^{(0)}, \cdots, [b_{j+1}^{(0)}, b_{j}^{(r)}] \cdots ]]. 
\]

\Pro \label{pro: PBW}
The algebra $\sTyar$ has a $\C$-basis consisting of ordered monomials in the elements
\[
\{ b_{\beta^\vee}^{(r)}\}_{\beta^\vee \in \Delta^+}^{r \geq 0} \cup 
\{z_{i}^{(2s_i+1)}\}_{i \in \indx}^{s_i \geq - \ai(\mu)}. 
\]
The same set is also a $\C[\hbar]$-basis of $\sTya$.
\enpro

\Proof
This is proven in \cite[Theorem 4.12]{LWZ1} for the unshifted case, and in \cite[Theorem 8.2, Proposition 8.9]{LPTTW} and \cite[Theorem 3.2]{TT} for the shifted case. For similar arguments in the untwisted case, see, e.g., \cite[Corollary 3.15]{Toda} and \cite[Theorem 2.55]{FT}. 
\enproof 

As a consequence of Theorem \ref{pro: PBW}, for any $\eta \in X_*(T)^+$ with $\eta \leq \mu$, the \emph{shift homomorphism} 
\eq \label{eq: shift homo}
{}^{\mathsf{tw}}\iota(\mu,\eta) \colon \sTya \to {}^{\mathsf{tw}}{{Y}}_{{\mu}-{\eta},\hbar}, \qquad 
b_{i}^{(r)} \mapsto b_{i}^{(r+\ai(\eta))}, \quad z_{i}^{(r)} \mapsto z_{i}^{(r+2\ai(\eta))} 
\eneq
is injective. If $\eta = \mu$, we abbreviate ${}^{\mathsf{tw}}\iota(\mu) = {}^{\mathsf{tw}}\iota(\mu,\mu)$. 
In particular, $\sTya$ can be realized as a subalgebra of the unshifted twisted Yangian $\sTy$ via ${}^{\mathsf{tw}}\iota(\mu)$. 


\subsection{Canonical filtration} 

We are primarily interested in the following subalgebra of the shifted twisted Yangian. 

\Defi
Let $\dya$ be the $\C[\h]$-subalgebra of $\sTya$ generated by 
\[
\{ \h b_{\beta^\vee}^{(r)}\}_{\beta^\vee \in \Delta^+}^{r \geq 0} \cup 
\{ \h z_{i}^{(2s_i+1)}\}_{i \in \indx}^{s_i \geq - \ai(\mu)}. 
\]
\enDefi 

We will now give an alternative description of $\dya$ as a Rees algebra, in analogy to the untwisted case. 
There is a filtration $F^\bullet_{\mu}$ on $\sTya$ determined by 
\[
\deg b_{\beta^\vee}^{(r)} = \beta^\vee(\mu) + r +1, \qquad \deg z_{i}^{(r)} = 2 \alpha_i^\vee(\mu) + r+1. 
\] 


\Thm
For any $\mu \in X_*(T)^+$, there is a canonical $\C[\h]$-algebra isomorphism
\[
\dya \cong \Rees^{F^\bullet_{\mu}} \sTyar. 
\]
\enthm 

\Proof
The proof relies on a comparison of Rees algebras with respect to the loop and canonical filtrations, and is formally the same as the proof of \cite[Theorem A.32]{FT}. 
\enproof

Passing to the semi-classical limit, there are canonical isomorphisms
\[
\dya / \h \dya \cong \Rees^{F^\bullet_{\mu}} \sTyar / \h \Rees^{F^\bullet_{\mu}} \sTyar \cong \gr^{F^\bullet_{\mu}} \sTyar. 
\]
Moreover, by \cite[Lemma 3.3]{TT}, $\gr^{F^\bullet_{\mu}} \sTyar$ is a commutative algebra. 
It is also naturally a Poisson algebra, with the Poisson bracket given by 
\[
\{ a, b \} = \hbar\mi [\hat{a}, \hat{b}] \quad \mbox{mod} \ \h\dya, 
\]
for any lifts $\hat{a}, \hat{b}$. 
Finally, we remark that the shift homomorphisms \eqref{eq: shift homo} induce monomorphisms of Poisson algebras
\[
{}^{\mathsf{tw}}\iota(\mu):\dya / \h \dya \hookrightarrow {}^{\mathsf{tw}} \mathbf{Y}_{\mu -\eta} / \h  {}^{\mathsf{tw}} \mathbf{Y}_{\mu -\eta}. 
\] 

\subsection{Reflection equation presentation}\label{sec-RTT}

The \emph{RTT twisted Yangian} $\tYrtt$ of split type associated to $\gln$ is the algebra generated by $s_{ij}^{(k)}$, where $k \geq 1$ and $1 \leq i,j \leq n$, subject to the relations 
\begin{align*}
(u^2 - v^2)[s_{ij}(u),s_{kl}(v)] &= \hbar (u+v) (s_{kj}(u)s_{il}(v) - s_{kj}(v)s_{il}(u)) \\ 
& - \hbar  (u-v) (s_{ik}(u)s_{jl}(v) - s_{ki}(v)s_{lj}(u)) \\ 
& + \hbar^2 (s_{ki}(u)s_{jl}(v) - s_{ki}(v)s_{jl}(u)) \\
s_{ji}(-u) &= s_{ij}(u) + \hbar \frac{s_{ij}(u) - s_{ij}(-u)}{2u}.  
\end{align*} 
Writing
\[
S(u) = \sum_{i,j=1}^n e_{ij} \otimes s_{ij}(u), \qquad 
s_{ij}(u) = \delta_{ij} + \sum_{k>0} s_{ij}^{(k)} u^{-k}
\]
these relations can also be presented in matrix form: 
\begin{align*}
R(u-v)S_1(u)R^t(-u-v)S_2(v) &= S_2(v) R^t(-u-v) S_1(u) R(u-v),  \\ 
S^t(-u) &= S(u) + \hbar \frac{S(u)-S(-u)}{2u}. 
\end{align*}

By \cite[Theorem 2.4.3]{Molev-book}, there is an injective algebra homomorphism
\[
\tYrtt \hookrightarrow \Yartt, \qquad S(u) \mapsto T^t(-u) T(u). 
\]
Moreover, by \cite[Theorem 2.10.1]{Molev-book}, this embedding endows $\tYrtt$ with the structure of a right coideal subalgebra. 

\Thm[{\cite[Corollary 2.4.4, Remark 2.4.5]{Molev-book}}]
Ordered monomials in the PBW variables 
\[
\{ s_{ij}^{(r)} \}_{i > j}^{r \geq 1} \cup \{ s_{ii}^{(2r)} \}_{i \in \indx}^{r \geq 1} 
\]
form a basis of $\tYrtt$ as a free $\C[\h]$-module, as well as a basis of ${}^{\mathsf{tw}} \widetilde{Y}^{\mathsf{rtt}}$ as a free $\C$-module. 
\enthm

We rely on the following key result linking the RTT and current realizations of the twisted Yangian. 

\Thm[{\cite[Theorem 5.1]{LWZ1}}] 
There is a canonical isomorphism 
\eq \label{eq: Upsilon} 
\Upsilon \colon  \Ty  \xrightarrow{\sim}  \tYrtt 
\eneq 
from the twisted Yangian $\Ty$ in the Drinfeld presentation (Definition \ref{def: tYangian}) to the RTT twisted Yangian $\tYrtt$, given by Gauss decomposition.  
\enthm 

Let ${}^{\mathsf{tw}} \widetilde{\mathbf{Y}}^{\mathsf{rtt}}$ be the $\C[\h]$-subalgebra of $\tYrtt$ generated by $\{ \h s_{ij}^{(r)} \}_{i,j \in \tilde{\indx}}^{r \geq 1}$. Recall that, when $\mu = 0$, we abbreviate ${}^{\mathsf{tw}} \widetilde{\mathbf{Y}} = {}^{\mathsf{tw}} \widetilde{\mathbf{Y}}_{0}$. 

\Pro \label{pro: Upsilon}
The isomorphism \eqref{eq: Upsilon} restricts to an isomorphism 
\[
\Upsilon \colon {}^{\mathsf{tw}} \widetilde{\mathbf{Y}}  \xrightarrow{\sim}  {}^{\mathsf{tw}} \widetilde{\mathbf{Y}}^{\mathsf{rtt}}. 
\]
\enpro

\Proof
The proof is entirely analogous to the proof of \cite[Proposition 2.34]{FT}. 
\enproof



\subsection{Ciccoli-Gavarini duality} 

When $\mu = 0$, the algebra $\dyazero = {}^{\mathsf{tw}} \mathbf{Y}_{0}$ can be identified with the Ciccoli--Gavarini dual of $\sTy$. 
Let us first recall the general definition of this concept. We continue to use the notation from \S \ref{sec: quantum duality}. 

Let $\mathfrak{b} \subset \mathfrak{a}$ be a Lie coideal, and suppose that $B$ is a \emph{compatible} deformation--quantization of $U(\mathfrak{b})$, i.e., $B \subset A$ is a coideal subalgebra over $\C[\h]$, and there is an isomorphism of Hopf algebras $B/\h B \cong U(\mathfrak{b})$. 
The \emph{Ciccoli--Gavarini dual} $B^\Lsh$ of $B$ is the following coideal subalgebra of $A'$: 
\[
B^\Lsh = \{ a \in B \mid \delta_m(a) \in \h^m A^{\otimes (m-1)} \otimes B  \mbox{ for all } m \in \N \} = B \cap A'. 
\] 
The importance of $B^\Lsh$ stems from the fact that, according to the quantum duality principle \cite[Theorem 3.3]{CG06}, $B^\Lsh$ is a deformation--quantization of the coordinate ring of the homogeneous space $G_{\mathfrak{a}^*}/G_{\mathfrak{b}^\perp}$.  

Let us now return to the specific setting of the twisted Yangian $\sTy$. It is a coideal subalgebra of $Y_\h$, and admits a grading with $\deg(\h) = 1$ and $\deg(x^{(r)}) = r$, for $x = z_i, b_{\beta^\vee}$ (or $\deg s_{ij}^{(r)} = r-1$ in the RTT presentation). 
The description of the coideal algebra structure can be found in, e.g., \cite[(2.64)]{Molev-book}. 
The twisted Yangian is a deformation--quantization of the universal enveloping algebra of the twisted current Lie algebra, i.e., $\sTy/ \h \sTy \cong U(\sln[t]^{\sigma})$, where $\sigma$ is the composition of the Cartan involution with the map $t \mapsto - t$. This isomorphism is an isomorphism of graded Hopf algebras if $U(\sln[t]^\sigma)$ is endowed with the loop grading. 

\Thm
There is a canonical $\C[\h]$-algebra isomorphism 
\[
\sTy^\Lsh \cong \dyazero. 
\]
\enthm

\Proof
The claim will follow from Proposition \ref{pro: Upsilon} if we can show ${}^{\mathsf{tw}}\widetilde{Y}_\h^\Lsh \cong {}^{\mathsf{tw}} \widetilde{\mathbf{Y}}^{\mathsf{rtt}}$. For the latter, one can apply the same argument as in \cite[Theorem A.26]{FT}. 
\enproof


\section{Sklyanin minors and ABCD presentation} 
\label{sec: Sklyanin}

Drinfeld's new realization of the Yangian admits a natural description in terms of quantum minors \cite{Drinfeld, NT, GKLO}. 
Below we deduce an analogous description of the current realization of the twisted Yangian from \cite{LWZ1} in terms of \emph{Sklyanin minors}. We assume $\mu = 0$ throughout.

\subsection{Current generators as Sklyanin minors}\label{sec-generators as minors}
For the sake of conciseness, we freely use the standard definition, properties and notation for the Sklyanin determinant and Sklyanin minors, as in, e.g., \cite[\S 2.5--\S2.6]{Molev-book}, without recalling these in detail. 

Recall the series $f_{i}(u), e_{i}(u), d_i(u)$ constructed by Gauss decomposition in \cite[(3.1)--(3.3)]{LWZ1}. 
The following proposition allows us to express the current generators of the twisted Yangian in terms of Sklyanin minors. 

\Pro
\label{lem: Sk2}
We have: 
\begin{align*}
f_{i}(u) &= s^{1\hdots i-1,i+1}_{1\hdots i}(u+\hbar(i-1)) (s^{1\hdots i}_{1 \hdots i}(u+\hbar(i-1)))\mi, \\
e_{i}(u) &= (s^{1\hdots i}_{1 \hdots i}(u+\hbar(i-1)))\mi s^{1\hdots i}_{1\hdots i-1,i+1}(u+\hbar(i-1)), \\ 
d_i(u) &= s^{1\hdots i}_{1 \hdots i}(u+\hbar(i-1))  (s^{1\hdots i-1}_{1 \hdots i-1}(u+\hbar(i-1)))\mi. 
\end{align*}
Hence
\begin{align*}
\textstyle b_i(u) &= f_{i}(u-\textstyle\frac{\hbar i}{2}) = s^{1\hdots i-1,i+1}_{1\hdots i}(u+\textstyle\frac{\hbar i}{2}-\h) (s^{1\hdots i}_{1 \hdots i}(u+\textstyle\frac{\hbar i}{2}-\h))\mi, \\
z_i(u) &= 
\frac{d_{i+1}(u-\textstyle\frac{\hbar i}{2})}{d_i(u-\textstyle\frac{\hbar i}{2})} 
= \frac{s^{1\hdots i-1}_{1 \hdots i-1}(u+\textstyle\frac{\hbar i}{2}-\h)  s^{1\hdots i+1}_{1 \hdots i+1}(u+\textstyle\frac{\hbar i}{2})}{s^{1\hdots i}_{1 \hdots i}(u+\textstyle\frac{\hbar i}{2}-\h)s^{1\hdots i}_{1 \hdots i}(u+\textstyle\frac{\hbar i}{2})}. 
\end{align*}
\enpro

\Proof
This follows immediately from the formulae in the proof of \cite[Corollary 3.2]{LWZ1} and \cite[Lemma 2.14.2]{Molev-book}. 
\enproof

We remark that, since $e_{i}(u) = f_{i}(-u-\h i)$, we also get $b_i(-u) = e_{i}(u-\frac{\h i}{2})$. 

\subsection{Relations} 

In this subsection, we will describe some of the relations between the Sklyanin minors. 
Define: 
\begin{align}
A_i(u) &= \textstyle s^{1\hdots i}_{1\hdots i}(u+\frac{\h i}{2}-\hbar), \quad& B_i(u) &= \textstyle s^{1\hdots i}_{1\hdots i-1,i+1}(u+\frac{\h i}{2}-\hbar), \\
C_i(u) &= \textstyle s^{1\hdots i-1,i+1}_{1\hdots i}(u+\frac{\h i}{2}-\h), \quad& D_i(u) &= \textstyle s^{1\hdots i-1,i+1}_{1\hdots i-1,i+1}(u+\frac{\h i}{2}-\h). 
\end{align}
Then, by Proposition \ref{lem: Sk2}, 
\[
B_i(u) = A_i(u) b_i(-u), \quad C_i(u) = b_i(u)A_i(u), \quad z_i(u) = 
\frac{A_{i+1}(u + \half \hbar) A_{i-1}(u+ \half \hbar)}{A_i(u) A_i(u + \h)}. 
\] 
Note that $A_i(u)$ satisfies $A_i(-u) = A_i(u+\h)$. 
If we set $\widetilde{A}_i(u) = A_i(u+\half \hbar)$, then 
$\widetilde{A}_i(-u) = \widetilde{A}_i(u)$ and 
\eq \label{eq: z in A terms}
z_i(u) = 
\frac{\widetilde{A}_{i+1}(u) \widetilde{A}_{i-1}(u)}{\widetilde{A}_i(u - \half \h) \widetilde{A}_i(u + \half\h)}. 
\eneq
We will similarly denote $\widetilde{X}_i(u) = X_i(u+ \half \h)$ for $X \in \{B,C,D\}$. 

Recall that the Sklyanin comatrix $\sco$ is defined by 
\[
\sco S(u-n+\h) = \sdet S(u). 
\]

\Lem
\label{lem: Sk1}
The following entries of the Sklyanin comatrix are equal to the corresponding Sklyanin minors: 
\begin{align}
\hs_{n-1,n-1}(u) &= s^{1\hdots n-2,n}_{1\hdots n-2,n}(u), \quad& \hs_{n,n}(u) &= s^{1\hdots n-1}_{1\hdots n-1}(u),  \\ 
\hs_{n-1,n}(u) &= s^{1\hdots n-2,n}_{1\hdots n-1}(u), \quad& \hs_{n,n-1}(u) &= s^{1\hdots n-1}_{1\hdots n-2,n}(u). 
\end{align}
\enlem 

\Proof
This is proven in the same way as \cite[Proposition 4.1]{Jing}. 
\enproof 

\Lem[{\cite[Proposition 2.12.3]{Molev-book}}]
\label{lem: Sk3}
The map 
\[
\textstyle
S(u) \mapsto \widehat{S}(-u+\frac{\h n}{2}-\h)
\]
is an automorphism. 
\enlem

Let $\dot{X}_i(u) = X_i(-u+ \half \h)$ for $X \in \{A,B,C,D\}$. 
Applying the two lemmas above, we get the following analogue of \cite[Proposition 2.1]{GKLO} and \cite[Proposition 1.2]{NT}. 
Note that the list below is not necessarily exhaustive, i.e., we expect more relations are needed to get a full presentation. 

\Pro
The following relations hold: 
\begin{align}
(u^2 - v^2)[\dot{A}_{i}(u), \dot{B}_{i}(v)] &= \hbar (u+v) (\dot{A}_i(u)\dot{B}_i(v) - \dot{A}_i(v)\dot{B}_i(u)) \\
&- \hbar (u-v) (\dot{A}_i(u)\dot{B}_i(v) - \dot{A}_i(v)\dot{C}_i(u)) \\ 
&+ \hbar^2 (\dot{A}_i(u)\dot{B}_i(v) - \dot{A}_i(v)\dot{B}_i(u)), \\ 
(u^2 - v^2)[\dot{A}_{i}(u), \dot{C}_{i}(v)] &= \hbar (u+v) (\dot{A}_i(u)\dot{C}_i(v) - \dot{A}_i(v)\dot{C}_i(u)) \\
&- \hbar (u-v) (\dot{A}_i(u)\dot{C}_i(v) - \dot{A}_i(v)\dot{B}_i(u)) \\ 
&+ \hbar^2 (\dot{A}_i(u)\dot{C}_i(v) - \dot{A}_i(v)\dot{C}_i(u)), \\ 
(u^2 - v^2)[\dot{B}_{i}(u), \dot{C}_{i}(v)] &= \hbar (u+v) (\dot{A}_i(u)\dot{D}_i(v) - \dot{A}_i(v)\dot{D}_i(u)) \\
&- \hbar (u-v) (\dot{B}_i(u)\dot{C}_i(v) - \dot{C}_i(v)\dot{B}_i(u)) \\ 
&- \hbar^2(u+v)\mi (\dot{A}_i(u)\dot{D}_i(v) - \dot{A}_i(v)\dot{D}_i(u)), \\ 
\dot{A}_{i}(-u) &= \dot{A}_{i}(u), \\
\dot{B}_i(-u) &= \dot{C}_i(u) + \hbar \frac{\dot{C}_i(u) - \dot{C}_i(-u)}{2u}, \\
\dot{C}_i(-u) &= \dot{B}_i(u) + \hbar \frac{\dot{B}_i(u) - \dot{B}_i(-u)}{2u}, \\
(u+v) [\dot{C}_i(u), \dot{C}_i(v)] &= \h \big(\dot{A}_i(v)\dot{D}_i(u) - \dot{A}_i(u)\dot{D}_i(v) \big), \\
(u+v) [\dot{B}_i(u), \dot{B}_i(v)] &= \h \big( \dot{D}_i(v)\dot{A}_i(u) - \dot{D}_i(u)\dot{A}_i(v) \big), \\ 
[\dot{A}_i(u), \dot{A}_i(v)] &= 0. 
\end{align}
\enpro

\Proof
This follows from Lemmas \ref{lem: Sk1}--\ref{lem: Sk3} and the RTT relations in the twisted Yangian (see, e.g., \cite[Proposition 2.2.1]{Molev-book}). 
\enproof


\section{Twisted GKLO representations} 
\label{sec: tGKLO}

In \cite{GKLO}, a remarkable family of representations of Yangians, quantizing moduli spaces of monopoles, was introduced by Gerasimov, Kharchev, Lebedev and Oblezin. Today they are commonly known as `GKLO representations'. 
They were later generalized to the cases of dominantly \cite{KWWY} and arbitrarily \cite{BFN} shifted Yangians. 
Below we construct analogous representations for the shifted \emph{twisted} Yangians $\Tya$.

\subsection{Difference operators} 

Fix a dominant coweight $\lambda \in X_*(T)^+$ with $\mu \leq \lambda$, and set $m_i = \omega_{i}^\vee(\lambda-\mu)$ and $\lambda_i = \alpha_{i}^\vee(\lambda)$. In other words, 
\[
\lambda = \sum_{i=1}^{n-1} \lambda_i \omega_{i}, \qquad 
\mu = \sum_{i=1}^{n-1} \mu_i \omega_{i}, \qquad 
\lambda - \mu = \sum_{i=1}^{n-1} m_i \alpha_{i}.  
\]  
Note that, since $\mu \leq \lambda$, each $m_i$ is a non-negative integer. 
Moreover, we have the identity
\eq \label{eq: lambda mu m}
\lambda_i - \mu_i = 2 m_i - \sum_{j \sim i}m_j, 
\eneq 
where $j \sim i$ means that $j$ and $i$ are connected by an edge in the Dynkin diagram. 
Let $\hdiff$ be the $\hcor$-algebra generated by
\begin{itemize} \setlength\itemsep{0.5em}
\item $\gamma_{i,k}, \beta_{i,k}^{\pm 1}$ (for $i \in \indx$ and $1 \leq k \leq m_i$), 
\item $((\gamma_{i,k} +r\h)^2 - (\gamma_{i,l}+s \h)^2)^{-1}$,  (for $k \neq l$ and $r,s \in \Z$), 
\item $(\gamma_{i,k} \pm \frac{r}{2} \hbar)^{-1}$ (for $r \in \Z$),
\end{itemize}
subject to the relations
\[
[\beta_{i,k}^{\pm1},\gamma_{j,l}] = \pm \delta_{ij}\delta_{kl}\hbar \beta_{i,k}^{\pm1}, \qquad [\gamma_{i,k}, \gamma_{j,l}] = 0 = [\beta_{i,k}, \beta_{j,l}], \qquad \beta_{i,k}^{\pm1}\beta_{i,k}^{\mp1} = 1. 
\]
The algebra $\hdiff$ has a natural representation on the space 
\[ \Pol = \hcor [\gamma_{i,k}, ((\gamma_{i,k} +r\h)^2 - (\gamma_{i,l}+s \h)^2)^{-1}, (\gamma_{i,k} + \textstyle\frac{r}{2} \hbar)^{-1}]_{i \in \indx, 1\leq k \neq l \leq m_i}^{r,s\in\Z},\]
where $\gamma_{i,k}$ acts by multiplication and $\beta_{i,k}^{\pm1}$ acts by the difference operator $e^{\pm \hbar \partial_{\gamma_{i,k}}}$. In analogy to \cite[Proposition 4.4]{KWWY}, there is an isomorphism of Poisson algebras
\[
\hdiff/\h\hdiff \cong 
\C [\gamma_{i,k}, \beta_{i,k}^{\pm1}, ((\gamma_{i,k} +r\h)^2 - (\gamma_{i,l}+s \h)^2)^{-1}, (\gamma_{i,k} \pm \textstyle\frac{r}{2} \hbar)^{-1}]_{i \in \indx, 1\leq k \neq l \leq m_i}^{r,s\in\Z}, 
\]
where the only non-trivial Poisson bracket between the generators on the RHS is $\{\beta_{i,k}^{\pm1}, \gamma_{i,k} \} = \pm \beta_{i,k}^{\pm1}$.

\subsection{Auxiliary relations} 

We begin by defining some auxiliary operators in $\hdiff$ and describing the relations between them. 
Let us abbreviate $\xi_{i,k} = \gamma_{i,k} + \hbar$. 
Choose polynomials 
\[
R_i(u) = \prod_{k=1}^{\lambda_i} (u^2 - r_{i,k}^2)
\qquad (i \in \indx),
\]
where $r_{i,k}$ are arbitrary complex numbers. 
Define 
\begin{align} \label{eq: xik defi}
\varkappa_{i,k} &= \frac{\prod_{l=1}^{m_{i+1}} (\gamma_{i,k}^2 - (\gamma_{i+1,l}+\frac12 \hbar)^2)\prod_{l=1}^{m_{i-1}} (\gamma_{i,k} + \gamma_{i-1,l}+\frac12 \hbar)}{2(\gamma_{i,k} - \half \hbar)\prod_{l \neq k}(\gamma_{i,k}^2 - \gamma_{i,l}^2)} \beta_{i,k}^{-1} \in \hdiff, \\ \label{eq: xik defi2}
\varkappa'_{i,k} &= R_i(\xi_{i,k}) \frac{\prod_{l=1}^{m_{i-1}} (\xi_{i,k} - (\gamma_{i-1,l}+\frac12 \hbar))}{2(\gamma_{i,k} + \frac32 \hbar)\prod_{l \neq k}(\xi_{i,k}^2 - \xi_{i,l}^2)} \beta_{i,k} \in \hdiff. 
\end{align}

The following lemma follows easily by direct calculation. 

\Lem \label{lem: five key rels}
The following relations hold: 
\begin{align} 
\label{eq: x-gamma}
[\varkappa_{i,k}, \gamma_{j,l}] &= - \hbar \delta_{ij}\delta_{kl} \varkappa_{i,k}, \qquad 
[\varkappa'_{i,k}, \gamma_{j,l}] =  \hbar \delta_{ij}\delta_{kl} \varkappa'_{i,k}, \\ 
 [\varkappa_{i,k}, \varkappa_{j,l}] &= \frac{\half a_{ij} \hbar}{\gamma_{i,k}-\gamma_{j,l}} [\varkappa_{j,l}, \varkappa_{i,k}]_+, \qquad
 [\varkappa'_{i,k}, \varkappa'_{j,l}] = [\varkappa'_{j,l}, \varkappa'_{i,k}]_+ \frac{-\half a_{ij}\hbar}{\gamma_{i,k}-\gamma_{j,l}}, \label{eq: xx} \\
\label{eq: xx'}
[\varkappa_{i,k}, \varkappa'_{j,l}] &= \frac{\half a_{ij} \hbar}{\gamma_{i,k}+\gamma_{j,l} + \hbar} [\varkappa_{i,k}, \varkappa'_{j,l}]_+ \qquad \big((i,k) \neq (j,l)\big), \\  
\varkappa_{i,k}\varkappa'_{i,k} &= R_i(\gamma_{i,k}) \frac{\prod_{j=i\pm1}\prod_{l=1}^{m_{j}} (\gamma_{i,k}^2 - (\gamma_{j,l}+\frac12 \hbar)^2)}{4(\gamma_{i,k} - \half \hbar)(\gamma_{i,k} + \half \hbar)\prod_{l \neq k}(\gamma_{i,k}^2 - \gamma_{i,l}^2)(\gamma_{i,k}^2 - \xi_{i,l}^2)}, \label{eq: xx3}
\\
\varkappa'_{i,k}\varkappa_{i,k} &= R_i(\xi_{i,k}) \frac{\prod_{j=i\pm1}\prod_{l=1}^{m_{j}} (\xi_{i,k}^2 - (\gamma_{j,l}+\frac12 \hbar)^2)}{4(\xi_{i,k} - \half \hbar)(\xi_{i,k} + \half \hbar)\prod_{l \neq k}(\xi_{i,k}^2 - \gamma_{i,l}^2)(\xi_{i,k}^2 - \xi_{i,l}^2)}.   \label{eq: xx4}
\\
\end{align} 
\enlem 
Relations \eqref{eq: xx}--\eqref{eq: xx'} can also be reformulated as
\begin{align}
(\gamma_{i,k} - \gamma_{j,l} - \half a_{ij} \hbar)\x_{i,k}\x_{j,l} &= (\gamma_{i,k} - \gamma_{j,l} + \half a_{ij} \hbar)\x_{j,l}\x_{i,k}, \label{eq: nxx1} \\ 
\x'_{i,k}\x'_{j,l} (\gamma_{i,k} - \gamma_{j,l} + \half a_{ij} \hbar) &= \x'_{j,l}\x'_{i,k}(\gamma_{i,k} - \gamma_{j,l} - \half a_{ij} \hbar), \label{eq: nxx2} \\ 
(\gamma_{i,k} + \gamma_{j,l} +(1- \half a_{ij}) \hbar)\x_{i,k}\x'_{j,l} &= (\gamma_{i,k} + \gamma_{j,l} + (1+\half a_{ij}) \hbar)\x'_{j,l}\x_{i,k}. \label{eq: nxx3}
\end{align}
For later convenience, for $i \neq j$, abbreviate 
\begin{align}
S_{i,k}^{j,l} =   \frac{4(\gamma_{i,k} - \half \hbar)(\gamma_{i,k} + \half \hbar)\x_{i,k} \x'_{i,k}}{\gamma_{i,k}^2 - (\gamma_{j,l}+\frac12 \hbar)^2}, \qquad 
T_{i,k}^{j,l} =   \frac{4(\xi_{i,k} - \half \hbar)(\xi_{i,k} + \half \hbar)\x'_{i,k} \x_{i,k}}{\xi_{i,k}^2 - (\xi_{j,l}+\frac12 \hbar)^2}. 
\end{align}
\noeqref{eq: nxx2}

\subsection{Twisted GKLO homomorphism} 

Below we define an analogue of the \emph{GKLO homomorphism} for shifted twisted Yangians. 

\Thm
\label{thm: GKLO}
The assignment 
\begin{align} \label{eq: bi action defi}
b_i(u) &\mapsto \sum_{k=1}^{m_i} \frac{1}{u-\gamma_{i,k}} \varkappa_{i,k} + \frac{1}{u+\gamma_{i,k}+ \hbar} \varkappa'_{i,k}, \\
\label{eq: h defi}
h_i(u) &\mapsto u^{-4m_1} \prod_{j=1}^{i-1} R_j(u-{\textstyle \frac{i-1-j}{2}\hbar}) \frac{\prod_{k=1}^{m_i} (u^2 - (\gamma_{i,k}+ \frac12 \hbar)^2)}{\prod_{l=1}^{m_{i-1}}(u-\gamma_{i-1,l})(u+\gamma_{i-1,l}+\hbar)}
\end{align}
defines a homomorphism $\Phi^\lambda_\mu \colon \Tya \to \hdiff$. 
\enthm

We define the \emph{$\lambda$-truncated $\mu$-shifted twisted Yangian} $\ttya$ to be the image of $\Tya$ under $\Phi^\lambda_\mu$. We refer to the pullback by $\Phi^\lambda_\mu$ of the natural representation $\Pol$ of $\hdiff$ as a \emph{twisted GKLO representation}. Restriction also yields representations of $\sTya$ and $\dya$. Let $\sttya$ (resp.\ $\dyatr$) denote the image of $\sTya$ (resp.\ $\dya$) under $\Phi^\lambda_\mu$. 

\subsection{Relations check} 
\label{sec: relations check}

We prove Theorem \ref{thm: GKLO} by explicitly checking that $\Phi^\lambda_\mu$ preserves the relations from Definition \ref{def: tYangian}. 
Relations \eqref{rel: hh} and \eqref{rel: bibj} are immediate. 

\subsubsection{Range of powers} 

We first check that the operators on the RHS of \eqref{eq: bi action defi}--\eqref{eq: h defi}, when expanded in $u\mi$, have the correct range of powers. The RHS of \eqref{eq: bi action defi} lies in $u\mi\hdiff[[u\mi]]$, matching $b_i(u)$. On the other hand, the RHS of \eqref{eq: h defi} lies in $u^{p}\hdiff[[u\mi]]$, where 
\[
p =  2 \left( \sum_{j=1}^{i-1} \lambda_j + m_i - m_{i-1} \right) - 4m_1 = -2m_1 + 2 \sum_{j=1}^i \mu_j = -2\ei(\mu), 
\]
by \eqref{eq: lambda mu m}, and the coefficient at the top power is $1$, which matches $h_i(u)$.

\subsubsection{Relation \eqref{rel: h1}}

It follows directly from \eqref{eq: h defi} that 
\begin{align} \label{eq: ziu}
z_i(u) = (h_i(u-\half \hbar))^{-1} h_{i+1}(u) \mapsto R_i(u) \frac{\prod_{j=i\pm1}\prod_{k=1}^{m_{j}} (u^2 - (\gamma_{j,k}+ \half \hbar)^2)}{\prod_{k=1}^{m_i} (u^2-\gamma_{i,k}^2)(u^2-\xi_{i,k}^2)}. 
\end{align} 
Since \eqref{eq: ziu} involves only even powers of $u$, \eqref{rel: h1} follows. 

\subsubsection{Relation \eqref{rel: hb}}

Given \eqref{eq: x-gamma}, it suffices to prove \eqref{rel: hb} in the case where $m_j = 1$, $\lambda_j=0$ and $m_{j\pm1} = 0$. Therefore, in this subsection we will omit the second subscript on $\gamma$ and $\varkappa$. 

First consider the case $i = j$. Let us calculate the LHS of \eqref{rel: hb}: 
\begin{align*}
h_i(u) b_i(v) &= \frac{(u^2 - (\gamma_i + \half \hbar)^2)}{(v-\gamma_i)} \varkappa_i + \frac{(u^2 - (\gamma_i + \half \hbar)^2)}{(v+\gamma_i+\hbar)} \varkappa'_i, \\
b_i(v) h_i(u) &= \frac{(u^2 - (\gamma_i - \half \hbar)^2)}{(v-\gamma_i)} \varkappa_i + \frac{(u^2 - (\gamma_i + \frac32 \hbar)^2)}{(v+\gamma_i+\hbar)} \varkappa'_i, \\
[h_i(u), b_i(v)] &= \frac{-2\hbar\gamma_i}{(v-\gamma_i)} \varkappa_i + \frac{2\hbar(\gamma_i+\hbar)}{(v+\gamma_i+\hbar)} \varkappa'_i. 
\end{align*}
Next, we compute the RHS: 
\begin{align*}
h_i(u)(b_i(v) - b_i(-u-\half \hbar)) &= \frac{(u^2 - (\gamma_i + \half \hbar)^2) (u+v + \half \hbar)}{(v-\gamma_i)(u +\gamma_i + \half \hbar)} \varkappa_i \\
&\quad +  \frac{(u^2 - (\gamma_i + \half \hbar)^2) (u+v + \half \hbar)}{(v+\gamma_i +\hbar)(u -\gamma_i - \half \hbar)} \varkappa'_i, \\ 
(b_i(v) - b_i(u+\half \hbar))h_i(u) &= \frac{(u^2 - (\gamma_i - \half \hbar)^2) (u-v + \half \hbar)}{(v-\gamma_i)(u - \gamma_i + \half \hbar)} \varkappa_i \\
&\quad + \frac{(u^2 - (\gamma_i + \frac32 \hbar)^2) (u-v + \half \hbar)}{(v+\gamma_i+\hbar)(u + \gamma_i + \frac32 \hbar)} \varkappa'_i. 
\end{align*} 
Hence
\begin{align*}
\frac{\hbar h_i(u)(b_i(v) - b_i(-u-\half \hbar))}{u+v+\frac12 \hbar} - \frac{\hbar(b_i(v) - b_i(u+\half \hbar))h_i(u)}{u-v+\frac12 \hbar} = \frac{-2\hbar\gamma_i}{(v-\gamma_i)} \varkappa_i + \frac{2\hbar(\gamma_i+\hbar)}{(v+\gamma_i+\hbar)} \varkappa'_i,
\end{align*}
as required. 

Now consider the case $i = j+1$. Let us calculate the LHS of \eqref{rel: hb}: 
\begin{align*}
h_{j+1}(u) b_j(v) &= \frac{1}{(v-\gamma_j)(u-\gamma_j)(u+\gamma_j+\hbar)} \varkappa_j + \frac{1}{(v+\gamma_j+\hbar)(u-\gamma_j)(u+\gamma_j+\hbar)} \varkappa'_j, \\
b_j(v) h_{j+1}(u) &= \frac{1}{(v-\gamma_j)(u-\gamma_j+\hbar)(u+\gamma_j)} \varkappa_j + \frac{1}{(v+\gamma_j+\hbar)(u-\gamma_j-\hbar)(u+\gamma_j+2\hbar)} \varkappa'_j, \\
[h_{j+1}(u), b_j(v)] &= \frac{2\hbar \gamma_j}{(v-\gamma_j)(u^2-\gamma_j^2)((u+\hbar)^2-\gamma_j^2)} \varkappa_j \\ 
&\quad - \frac{2\hbar (\gamma_j + \hbar)}{(v+\gamma_j-\hbar)(u^2-(\gamma_j+\hbar)^2)(u-\gamma_j)(u+\gamma_j + 2\hbar)} \varkappa'_j. 
\end{align*} 
Next, we compute the RHS: 
\begin{align*}
h_j(u)\big(b_j(v) - b_j(-u)\big) &= \frac{u+v}{(v-\gamma_j)(u^2 - \gamma_j^2)(u+\gamma_j+\hbar)} \varkappa_j \\
&\quad +  \frac{u+v}{(v+\gamma_j +\hbar)(u-\gamma_j)(u^2-(\gamma_j+\hbar)^2)} \varkappa'_j, \\ 
\big(b_j(v) - b_j(u)\big)h_j(u) &= \frac{u-v}{(v-\gamma_j)(u^2-\gamma_j^2)(u-\gamma_j+\hbar)} \varkappa_j \\
&\quad +  \frac{u-v}{(v+\gamma_j+\hbar)(u^2-(\gamma_j+\hbar)^2)(u+\gamma_j+2\hbar)} \varkappa'_j. 
\end{align*}
Hence
\begin{align*}
&\frac{\hbar(b_j(v) - b_j(u))h_j(u)}{u-v} - \frac{\hbar h_j(u)(b_j(v) - b_j(-u))}{u+v} = \\
=&\frac{2\hbar \gamma_j}{(v-\gamma_j)(u^2-\gamma_j^2)((u+\hbar)^2-\gamma_j^2)} \varkappa_j
- \frac{2 \hbar (\gamma_j + \hbar)}{(v+\gamma_j-\hbar)(u^2-(\gamma_j+\hbar)^2)(u-\gamma_j)(u+\gamma_j + 2\hbar)} \varkappa'_j, 
\end{align*}
as required. 

\subsubsection{Relation \eqref{eq: bbz}} 
\label{sec: 15 rel}

We first calculate the LHS: 
\begin{align*}
[b_i(u),b_i(v)] &= \sum_k \frac{\hbar(v-u)}{(u-\gamma_{i,k})(v-\gamma_{i,k})(u-\gamma_{i,k}+\hbar)(v-\gamma_{i,k}+\hbar)} \varkappa_{i,k}^2 \\ 
&\quad +  \sum_k \frac{\hbar(v-u)}{(u+\gamma_{i,k}+\hbar)(v+\gamma_{i,k}+\hbar)(u+\gamma_{i,k}+2\hbar)(v+\gamma_{i,k}+2\hbar)} (\varkappa'_{i,k})^2 \\
&\quad + \sum_{k\neq l} \frac{(v-u)(\gamma_{i,k}-\gamma_{i,l})}{(u-\gamma_{i,k})(v-\gamma_{i,l})(v-\gamma_{i,k})(u-\gamma_{i,l})} \varkappa_{i,k} \varkappa_{i,l} \\ 
&\quad - \sum_{k \neq l} \varkappa'_{i,k} \varkappa'_{i,l} \frac{(v-u)(\gamma_{i,k}-\gamma_{i,l})}{(u+\gamma_{i,k})(v+\gamma_{i,l})(v+\gamma_{i,k})(u+\gamma_{i,l})} \\
&\quad + \sum_{k\neq l} \frac{(v-u)(\gamma_{i,k}+\gamma_{i,l}+\hbar)}{(u-\gamma_{i,k})(v-\gamma_{i,k})(u+\gamma_{i,l}+\hbar)(v+\gamma_{i,l}+\hbar)} [\varkappa_{i,k}, \varkappa'_{i,l}] \\ 
&\quad +  \sum_k \frac{2\gamma_{i,k}(v-u)}{(u^2-\gamma_{i,k}^2)(v^2-\gamma_{i,k}^2)} \varkappa_{i,k} \varkappa'_{i,k} \\
&\quad -  \sum_k \frac{2(\gamma_{i,k}+\hbar)(v-u)}{(u^2 - (\gamma_{i,k}+\hbar)^2)(v^2 - (\gamma_{i,k}+\hbar)^2)} \varkappa'_{i,k} \varkappa_{i,k}. 
\end{align*}
Let us denote each of the seven summands above as $X_1, \cdots, X_7$, counting from the top. 

We now pass to the RHS of \eqref{eq: bbz}: 
\begin{align*}
b_i(u) - b_i(v) &= \sum_k \frac{v-u}{(u-\gamma_{i,k})(v-\gamma_{i,k})}\varkappa_{i,k} 
+ \frac{v-u}{(u+\gamma_{i,k}+\hbar)(v+\gamma_{i,k}+\hbar)}\varkappa'_{i,k}, \\ 
(b_i(u) - b_i(v))^2 &=  \sum_k 
\frac{(v-u)^2}{(u-\gamma_{i,k})(v-\gamma_{i,k})(u-\gamma_{i,k}+\hbar)(v-\gamma_{i,k}+\hbar)}\varkappa_{i,k}^2 \\
&\quad + \sum_k \frac{(v-u)^2}{(u+\gamma_{i,k}+\hbar)(v+\gamma_{i,k}+\hbar)(u+\gamma_{i,k}+2\hbar)(v+\gamma_{i,k}+2\hbar)}(\varkappa'_{i,k})^2 \\ 
 &\quad +  \sum_{k,l} 
\frac{(v-u)^2}{(u-\gamma_{i,k})(v-\gamma_{i,k})(u-\gamma_{i,l})(v-\gamma_{i,l})}\varkappa_{i,k}\varkappa_{i,l} \\
 &\quad +  \sum_{k,l} \varkappa'_{i,k}\varkappa'_{i,l}
\frac{(v-u)^2}{(u+\gamma_{i,k})(v+\gamma_{i,k})(u+\gamma_{i,k})(v+\gamma_{i,k})} \\
&\quad + \sum_{k\neq l} \frac{(v-u)^2}{(u-\gamma_{i,k})(v-\gamma_{i,k})(u+\gamma_{i,l}+\hbar)(v+\gamma_{i,l}+\hbar)} [\varkappa_{i,k}, \varkappa'_{i,l}]_+ \\ 
&\quad + \sum_k \frac{(v-u)^2}{(u^2-\gamma_{i,k}^2)(v^2-\gamma_{i,k}^2)} \varkappa_{i,k} \varkappa'_{i,k} \\ 
&\quad + \sum_k \frac{(v-u)^2}{(u^2-(\gamma_{i,k}+\hbar)^2)(v^2-(\gamma_{i,k}+\hbar)^2)} \varkappa'_{i,k} \varkappa_{i,k}. 
\end{align*}  
Let us denote each of the seven summands above as $Y_1, \cdots, Y_7$, counting from the top. We see immediately that 
$X_1 = \frac{\hbar}{v-u}Y_1$ and $X_2 = \frac{\hbar}{v-u}Y_2$. 
Relations \eqref{eq: xx}--\eqref{eq: xx'} also imply that $X_3 = \frac{\hbar}{v-u}Y_3$, $X_4 = \frac{\hbar}{v-u}Y_4$, and  $X_5 = \frac{\hbar}{v-u}Y_5$.  
Therefore, to establish \eqref{eq: bbz}, it suffices to show that 
\begin{align} 
\label{eq: toshow}
\sum_k \frac{2(\gamma_{i,k}-\half \hbar)(v^2-u^2)}{(u^2-\gamma_{i,k}^2)(v^2-\gamma_{i,k}^2)} \varkappa_{i,k} \varkappa'_{i,k} 
&- \sum_k \frac{2(\gamma_{i,k}+\frac32\hbar)(v^2-u^2)}{(u^2 - \xi_{i,k}^2)(v^2 - \xi_{i,k}^2)} \varkappa'_{i,k} \varkappa_{i,k} = \\
&= \hbar(\uz_i(v) - \uz_i(-u)). 
\end{align}

Let us first calculate the LHS of \eqref{eq: toshow}. Using relations \eqref{eq: xx3}--\eqref{eq: xx4}, we get 
\begin{align} \label{eq: LHS}
\ \ \ \LHS &= \sum_k \left( \frac{1}{u^2 - \gamma_{i,k}^2}
- \frac{1}{v^2 - \gamma_{i,k}^2} \right) \frac{ R_i(\gamma_{i,k}) \prod_{j=i\pm1}\prod_{l=1}^{m_{j}} (\gamma_{i,k}^2 - (\gamma_{j,l}+\frac12 \hbar)^2)}{2(\gamma_{i,k} + \half \hbar)\prod_{l \neq k}(\gamma_{i,k}^2 - \gamma_{i,l}^2)(\gamma_{i,k}^2 - \xi_{i,l}^2)} \\
&+ \sum_k \left( \frac{1}{v^2 - \xi_{i,k}^2}
- \frac{1}{u^2 - \xi_{i,k}^2} \right) 
 \frac{R_i(\xi_{i,k})\prod_{j=i\pm1}\prod_{l=1}^{m_{j}} (\xi_{i,k}^2 - (\gamma_{j,l}+\frac12 \hbar)^2)}{2(\gamma_{i,k} + \half \hbar)\prod_{l \neq k}(\xi_{i,k}^2 - \gamma_{i,l}^2)(\xi_{i,k}^2 - \xi_{i,l}^2)}. 
\end{align}
Next, we compute the RHS of \eqref{eq: toshow}. 
Applying partial fraction decomposition to the denominator of \eqref{eq: ziu}, we get 
\begin{align} \label{eq: ziu2}
-\hbar z_i(u) =& \ \sum_{k=1}^{m_i}  \frac{R_i(u) \prod_{j=i\pm1}\prod_{l=1}^{m_{j}} (u^2 - (\gamma_{j,l} + \half \hbar)^2)}{2(\gamma_{i,k}+\half \hbar) \prod_{l \neq k} (\gamma_{i,k}^2-\gamma_{i,l}^2)(\gamma_{i,k}^2-\xi_{i,l}^2)}\frac{1}{u^2 - \gamma_{i,k}^2} \\ 
&\ -  \sum_{k=1}^{m_i}  \frac{R_i(u) \prod_{j=i\pm1}\prod_{l=1}^{m_{j}} (u^2 - (\gamma_{j,l} + \half \hbar)^2)}{2(\gamma_{i,k}+\half \hbar) \prod_{l \neq k} (\xi_{i,k}^2-\gamma_{i,l}^2)(\xi_{i,k}^2-\xi_{i,l}^2)} \frac{1}{u^2 - \xi_{i,k}^2}. 
\end{align} 
The result now follows by comparing the `$u$-part' of \eqref{eq: LHS} with the principal part of \eqref{eq: ziu2}, considered as Laurent series in $u^{-2}$ (with analogous comparison for the `$v$-part'). This follows from the fact that the principal parts of $\frac{p(\gamma_{i,k}^2)}{u^2-\gamma_{i,k}^2}$ and $\frac{p(u^2)}{u^2-\gamma_{i,k}^2}$ coincide for any polynomial $p(\cdot)$. 

\subsubsection{Relation \eqref{rel: bibi1}}

It suffices to prove \eqref{rel: bibi1} for $m_i = m_{i+1} = 1$, so we drop the second subscript on $\gamma$ and $\varkappa$. 
 The LHS is: 
\begin{align*}
 (u-v)[b_{i}(u), b_{i+1}(v)] &= \frac{u-v}{(u-\gamma_i)(v-\gamma_{i+1})} [\varkappa_i, \varkappa_{i+1}] 
 + [\varkappa'_i, \varkappa'_{i+1}] \frac{u-v}{(u+\gamma_i)(v+\gamma_{i+1})} \\ 
 &\, + \frac{u-v}{(u-\gamma_i)(v+\gamma_{i+1}+\hbar)}[\varkappa_i,\varkappa'_{i+1}] \\ 
&\, +  \frac{u-v}{(v-\gamma_{i+1})(u+\gamma_{i}+\hbar)}[\varkappa'_i,\varkappa_{i+1}].  
\end{align*}
On the other hand, on the RHS, 
\begin{align*}
-\half \hbar [b_i(u), b_{i+1}(v)]_+ &= 
\frac{-\half \hbar}{(u-\gamma_i)(v-\gamma_{i+1})} [\varkappa_i, \varkappa_{i+1}]_+
 + [\varkappa'_i, \varkappa'_{i+1}]_+ \frac{-\half \hbar}{(u+\gamma_i)(v+\gamma_{i+1})} \\ 
 &\, + \frac{-\half \hbar}{(u-\gamma_i)(v+\gamma_{i+1}+\hbar)}[\varkappa_i,\varkappa'_{i+1}]_+ \\
&\, +  \frac{-\half \hbar}{(v-\gamma_{i+1})(u+\gamma_{i}+\hbar)}[\varkappa'_i,\varkappa_{i+1}]_+, 
\end{align*}
and 
\begin{align*}
[b_{i}^{(0)}, b_{i+1}(v)]+[b_{i+1}^{(0)},b_i(u)] &= 
\frac{u-v+\gamma_{i+1} - \gamma_i}{(u-\gamma_i)(v-\gamma_{i+1})} [\varkappa_i, \varkappa_{i+1}] 
 + [\varkappa'_i, \varkappa'_{i+1}] \frac{u-v-\gamma_{i+1} + \gamma_i}{(u+\gamma_i)(v+\gamma_{i+1})} \\ 
 &\, + \frac{u-v-\gamma_i-\gamma_{i+1}-\hbar}{(u-\gamma_i)(v+\gamma_{i+1}+\hbar)}[\varkappa_i,\varkappa'_{i+1}] \\
&\, +  \frac{u-v+\gamma_i+\gamma_{i+1}+\hbar}{(v-\gamma_{i+1})(u+\gamma_{i}+\hbar)}[\varkappa'_i,\varkappa_{i+1}].  
\end{align*} 
The result is now implied by the identities
\begin{align*}
(\gamma_{i+1} - \gamma_i)[\varkappa_i, \varkappa_{i+1}] &= \half \hbar [\varkappa_i, \varkappa_{i+1}]_+, \quad & 
[\varkappa'_i, \varkappa'_{i+1}] (\gamma_{i} - \gamma_{i+1}) &= \half \hbar [\varkappa'_i, \varkappa'_{i+1}]_+, \\
(\gamma_{i+1} + \gamma_i + \hbar)[\varkappa_i, \varkappa'_{i+1}] &= -\half \hbar [\varkappa_i, \varkappa'_{i+1}]_+, \quad & 
(\gamma_{i+1} + \gamma_i + \hbar)[\varkappa'_i, \varkappa_{i+1}] &= \half \hbar [\varkappa_i, \varkappa'_{i+1}]_+,  
\end{align*}
which follow from \eqref{eq: xx}--\eqref{eq: xx'}. 

\subsubsection{Relation \eqref{rel: Serre}} 

Without loss of generality, we may assume $j=i+1$. Consider the LHS of \eqref{rel: Serre} as a (degree $3$) noncommutative polynomial in $\varkappa_{i,k}, \varkappa'_{i,l}, \varkappa_{i+1,m}, \varkappa'_{i+1,n}$. 
First, we show that the sum of monomials \emph{not} containing both $\varkappa_{i,k}$ and $\varkappa'_{i,k}$ (with the same second index) vanishes. Such monomials come in four different types: 
\begin{enumerate}[itemsep=2mm]
\item monomials containing $\x_{i,k}, \x_{i,l}, \x_{i+1,m}$ or $\x'_{i,k}, \x'_{i,l}, \x'_{i+1,m}$ (if $k=l$ then the associated variable occurs twice); 
\item monomials containing $\x_{i,k}, \x_{i,l}, \x'_{i+1,m}$ or $\x'_{i,k}, \x'_{i,l}, \x_{i+1,m}$ ($k\neq l$); 
\item monomials containing $\x_{i,k}$ with multiplicity $2$ and $\x'_{i+1,m}$, or $\x'_{i,k}$ with multiplicity $2$ and $\x_{i+1,m}$; 
\item monomials containing $\x_{i,k}, \x'_{i,l}$ ($k \neq l$) and $\x_{i+1,m}$ or $\x'_{i+1,m}$. 
\end{enumerate}
The vanishing of the sum of all monomials of type $1$ follows directly from the argument in the non-twisted case, i.e.,  \cite[Lemma 3.1]{GKLO} or \cite[B(vi)]{BFN}. For the other cases, 
we will need the following lemma. 
 
\Lem \label{lem: mixed 3}
The following identities hold: 
\begin{align*}
[\x_{i,k},[\x_{i,k}, \x'_{i+1,m}]] &= 0, \\ 
[\x_{i,k},[\x_{i,l}, \x'_{i+1,m}]] &= \frac{-\h^2(\gamma_{i,k}+\gamma_{i,l}+2\gamma_{i+1,m} + 2\h)}{(\gamma_{i,k}+\gamma_{j,m}+ \frac{\h}{2})(\gamma_{i,l}+\gamma_{j,m}+\frac{\h}{2})(\gamma_{i,k}-\gamma_{i,l}+\h)} \x_{i,k} \x_{i,l} \x'_{i+1,m}, \\ 
[\x_{i,k}, [\x'_{i,l}, \x_{i+1,m}]] &= 
\frac{\h^2(\gamma_{i,l} - \gamma_{i,k} + 2\gamma_{i+1,m} + \h)}{(\gamma_{i+1,m}-\gamma_{i,k}+\frac{\h}{2})(\gamma_{i,l}+\gamma_{i,k}+2\h)(\gamma_{i+1,m}+\gamma_{i,l}+\frac{3\h}{2})} \x_{i,k}\x'_{i,l}\x_{i+1,m}, \\
[\x'_{i,k}, [\x_{i,l}, \x_{i+1,m}]] &= 
\frac{-\h^2(\gamma_{i,k} - \gamma_{i,l} + 2\gamma_{i+1,m} + \h)}{(\gamma_{i+1,m}-\gamma_{i,l}+\frac{\h}{2})(\gamma_{i,l}+\gamma_{i,k})(\gamma_{i+1,m}+\gamma_{i,k}+\frac{3\h}{2})} \x'_{i,k}\x_{i,l}\x_{i+1,m}, 
\end{align*}
for $k \neq l$. 
\enlem 

\begin{proof}
The lemma follows by direct calculation using \eqref{eq: nxx1}--\eqref{eq: nxx3}. 
\end{proof} 

Since 
\begin{align} \label{eq: br reduct}
&\Sym_{u_1, u_2} \left[ \frac{1}{u_1-\gamma_{i,k}}\x_{i,k}, \left[  \frac{1}{u_2-\gamma_{i,k}}\x_{i,k}, \frac{1}{t+\gamma_{i+1,m}+\h}\x'_{i+1,m}\right]\right] = \\
&= \Sym_{u_1, u_2}  \frac{1}{(u_1-\gamma_{i,k}+\h)(u_2 - \gamma_{i,k})(t+\gamma_{i,k}+\h)} \big[\x_{i,k},\big[\x_{i,k}, \x'_{i+1,m}\big]\big], 
\end{align} 
the first formula of Lemma \ref{lem: mixed 3} implies that the sum of all monomials of type $3$, containing 
$\x_{i,k}$ with multiplicity $2$ and $\x'_{i+1,m}$, 
 on the LHS of \eqref{rel: Serre}, vanishes. 
 It is clear that the other subcase, i.e., monomials of type $3$ containing $\x'_{i,k}$ with multiplicity $2$ and $\x_{i+1,m}$, can be handled using an analogous argument. 

In the other cases, a similar calculation to \eqref{eq: br reduct} also shows that one can ignore the denominators, such as $\frac{1}{u_1 - \gamma_{p,s}}$, in front of $\x_{p,s}, \x'_{r,t}$.  
Let $k \neq l$. Then 
$
\x_{i,l} \x_{i,k} \x'_{i+1,m} = \frac{\gamma_{i,k}-\gamma_{i,l}-\h}{\gamma_{i,k}-\gamma_{i,l}+\h} \x_{i,k} \x_{i,l} \x'_{i+1,m}, 
$
and the second formula of Lemma \ref{lem: mixed 3} implies that \linebreak $[\x_{i,k},[\x_{i,l}, \x'_{i+1,m}]] + [\x_{i,l},[\x_{i,k}, \x'_{i+1,m}]] = 0$. 
Hence the sum of all monomials of type $2$ on the LHS of \eqref{rel: Serre} vanishes. 

Again, let $k \neq l$. Then 
$\x'_{i,l}\x_{i,k}\x_{i+1,m} = \frac{\gamma_{i,k}+\gamma_{i,l}}{\gamma_{i,k}+\gamma_{i,l}+2\h} \x_{i,k}\x'_{i,l}\x_{i+1,m}$, and the third and fourth formulae of Lemma \ref{lem: mixed 3} imply that 
$[\x_{i,k}, [\x'_{i,l}, \x_{i+1,m}]] + [\x'_{i,l}, [\x_{i,k}, \x_{i+1,m}]] =0$. It follows that the sum of all monomials of type $4$ on the LHS of \eqref{rel: Serre} also vanishes.

We will now prove that the sum of the remaining monomials equals the RHS of \eqref{rel: Serre}. We will need the following lemma. 
 
\Lem \label{lem: xxx ST}
The following identities hold:  
\begin{align} 
\x_{i,k} \x'_{i,k} \x_{i+1,l} - 2 \x_{i,k} \x_{i+1,l} \x'_{i,k} +  \x_{i+1,l} \x_{i,k} \x'_{i,k} &=  \frac{\half \hbar S_{i,k}^{i+1,l}}{\gamma_{i,k} + \half \hbar} \varkappa_{i+1,l},\\ 
\x'_{i,k} \x_{i,k} \x_{i+1,l} - 2 \x'_{i,k} \x_{i+1,l} \x_{i,k} +  \x_{i+1,l} \x'_{i,k} \x_{i,k} &=  \frac{\half \hbar T_{i,k}^{i+1,l}}{\gamma_{i,k} + \half \hbar} \varkappa_{i+1,l}, \\ 
\x_{i,k} \x'_{i,k} \x'_{i+1,l} - 2 \x_{i,k} \x'_{i+1,l} \x'_{i,k} +  \x'_{i+1,l} \x_{i,k} \x'_{i,k} &=  \frac{\half \hbar S_{i,k}^{i+1,l}}{\gamma_{i,k} + \half \hbar} \varkappa'_{i+1,l}, \\ 
\x'_{i,k} \x_{i,k} \x'_{i+1,l} - 2 \x'_{i,k} \x'_{i+1,l} \x_{i,k} +  \x'_{i+1,l} \x'_{i,k} \x_{i,k} &=  \frac{\half \hbar T_{i,k}^{i+1,l}}{\gamma_{i,k} + \half \hbar} \varkappa'_{i+1,l}. \\ 
\end{align}
\enlem

\Proof
We have
\begin{align}
\x_{i,k} \x'_{i,k} \x_{i+1,l} &= 
\frac{\gamma_{i,k}^2 - (\gamma_{i+1,l}+\frac12 \hbar)^2}{4(\gamma_{i,k} - \half \hbar)(\gamma_{i,k} + \half \hbar)} S_{i,k}^{i+1,l}  \x_{i+1,l}, \\
\x_{i,k} \x_{i+1,l} \x'_{i,k}  &= 
\frac{(\gamma_{i,k} -\frac12 \hbar)^2 - \gamma_{i+1,l}^2 }{4(\gamma_{i,k} - \half \hbar)(\gamma_{i,k} + \half \hbar)} S_{i,k}^{i+1,l}  \x_{i+1,l}, \\
\x_{i+1,l} \x_{i,k} \x'_{i,k}  &= 
\frac{\gamma_{i,k}^2 - (\gamma_{i+1,l}-\frac12 \hbar)^2}{4(\gamma_{i,k} - \half \hbar)(\gamma_{i,k} + \half \hbar)} S_{i,k}^{i+1,l}  \x_{i+1,l}.
\end{align}
The first identity of the lemma now follows from the fact that
\[
\big(\gamma_{i,k}^2 - (\gamma_{i+1,l}+\half \hbar)^2\big) - 2\big((\gamma_{i,k} -\half \hbar)^2 - \gamma_{i+1,l}^2 \big) + \big( \gamma_{i,k}^2 - (\gamma_{i+1,l}-\half \hbar)^2\big) 
= 2\hbar (\gamma_{i,k} - \half \hbar).  
\]
An analogous argument establishes the other identities. 
\enproof 

Lemma \ref{lem: xxx ST} implies that the LHS of \eqref{rel: Serre} is equal to
\[
\LHS = \sum_{u = u_1,u_2} \left( \sum_{k,l} \frac{\hbar u S_{i,k}^{i+1,l}}{(u^2-\gamma_{i,k}^2)(\gamma_{i,k} + \half \hbar)} b_{i+1}^l(t) 
+  \sum_{k,l} \frac{\hbar u T_{i,k}^{i+1,l}}{(u^2-\xi_{i,k}^2)(\gamma_{i,k} + \half \hbar)} b_{i+1}^l(t) \right), 
\]
where
\[
b_{i+1}^l(t) = \frac{1}{t-\gamma_{i+1,l}} \varkappa_{i+1,l} + \frac{1}{t+\gamma_{i+1,l}+ \hbar} \varkappa'_{i+1,l}.
\]
On the other hand, 
using \eqref{eq: ziu}, and the fact that 
\begin{align}
&\left( (t-\h )\, z_i(u) \x_{i+1,l}- (t+\h )\, \x_{i+1,l}z_i(u) \right) = \\ 
=& \frac{\hbar\big((-2u^2+\half \h^2) + 2\gamma_{i+1,l}(\gamma_{i+1,l}-t)\big)}{u^2 - (\gamma_{i+1,l} + \half \hbar)^2} z_i(u)\x_{i+1,l}, \\ 
&\left( (t-\h )\, z_i(u) \x'_{i+1,l}- (t+\h )\, \x'_{i+1,l}z_i(u) \right) = \\ 
=& \frac{\hbar\big((-2u^2+\half \h^2) + 2\xi_{i+1,l}(\xi_{i+1,l}+t)\big)}{u^2 - (\gamma_{i+1,l} + \half \hbar)^2} z_i(u)\x'_{i+1,l}, 
\end{align}
we conclude that 
\begin{align}
 &\frac{(t-\h )\, z_i(u) b_{i+1}(t)- (t+\h )\, b_{i+1}(t)z_i(u) }{4u^2-\h^2} 
 =_t \\ 
 &=_t - \sum_{l} R_i(u) \frac{\prod_{p=1}^{m_{i-1}} (u^2 - (\gamma_{i-1,p}+ \half \hbar)^2)\prod_{l \neq p=1}^{m_{i+1}} (u^2 - (\gamma_{i+1,k}+ \half \hbar)^2)}{2\prod_{k=1}^{m_i} (u^2-\gamma_{i,k}^2)(u^2-\xi_{i,k}^2)} b_{i+1}^l(t), 
\end{align}  
where $=_t$ denotes the equality of $t$-principal parts. 
A partial fraction decomposition argument, analogous to \S \ref{sec: 15 rel}, completes the proof of \eqref{rel: Serre}. 
This also concludes the proof of Theorem \ref{thm: GKLO}. 

\subsection{Symmetrized representation}

Our definition of the operators \eqref{eq: xik defi}--\eqref{eq: xik defi2} is obviously asymmetric. This asymmetry can be resolved at the cost of working in a certain quadratic extension of $\hdiff$. 
More precisely, let 
$\hdiffs$ be the $\hcor$-algebra generated by
\begin{itemize} \setlength\itemsep{0.5em}
\item $\gamma_{i,k}, \beta_{i,k}^{\pm 1}$ (for $i \in \indx$ and $1 \leq k \leq m_i$), 
\item $((\gamma_{i,k} +r\h)^2 - (\gamma_{i,l}+s \h)^2)^{-1}$,  (for $k \neq l$ and $r,s \in \Z$), 
\item $(\gamma_{i,k} \pm \frac{r}{2} \hbar)^{-1}$ (for $r \in \Z$), 
\item $(\gamma_{i,k} - \gamma_{i \pm 1,l} +\frac{r}{2}\h)^{\frac12}$ and $(\gamma_{i,k} + \gamma_{i \pm 1,l} +\frac{r}{2}\h)^{\frac12}$,  (for $r \in \Z$), 
\end{itemize}
subject to the relations
\[
[\beta_{i,k}^{\pm1},\gamma_{j,l}] = \pm \delta_{ij}\delta_{kl}\hbar \beta_{i,k}^{\pm1}, \qquad [\gamma_{i,k}, \gamma_{j,l}] = 0 = [\beta_{i,k}, \beta_{j,l}], \qquad \beta_{i,k}^{\pm1}\beta_{i,k}^{\mp1} = 1. 
\] 
Let us fix a root $\sqrt{-1}$ of $-1$. In the definition above, 
we choose the roots consistently, so that, e.g., 
\[
(\gamma_{i,k} - \gamma_{i+1,l})^{\frac12} = \sqrt{-1} (\gamma_{i+1,l} - \gamma_{i,k})^{\frac12}. 
\]
Then we can define 
\begin{align} 
\breve{\varkappa}_{i,k} &= \left(R_i(\gamma_{i,k})\right)^{\frac12} \frac{\prod_{j \in \{i\pm1\}}\prod_{l=1}^{m_{i-1}} (\gamma_{i,k}^2 - (\gamma_{j,l}+\frac12 \hbar)^2)^{\frac12}}{2(\gamma_{i,k} - \half \hbar)\prod_{l \neq k}(\gamma_{i,k}^2 - \gamma_{i,l}^2)} \beta_{i,k}^{-1} \in \hdiffs, \\ 
\breve{\varkappa}'_{i,k} &= \left(R_i(\xi_{i,k})\right)^{\frac12} \frac{\prod_{j \in \{i\pm1\}}\prod_{l=1}^{m_{i-1}} (\xi_{i,k}^2 - (\gamma_{j,l}+\frac12 \hbar)^2)^{\frac12}}{2(\gamma_{i,k} + \frac32 \hbar)\prod_{l \neq k}(\xi_{i,k}^2 - \xi_{i,l}^2)} \beta_{i,k} \in \hdiffs. 
\end{align} 
One easily checks that Lemma \ref{lem: five key rels} still holds if we replace $\x_{i,k} \leftrightarrow \breve{\x}_{i,k}$ and $\x'_{i,k} \leftrightarrow \breve{\x}'_{i,k}$. In particular, carrying out appropriate modifications throughout \S \ref{sec: relations check}, we get the following version of Theorem \ref{thm: GKLO}. 

\Cor
\label{thm: GKLO sym}
The assignment 
\begin{align}
b_i(u) &\mapsto \sum_{k=1}^{m_i} \frac{1}{u-\gamma_{i,k}} \breve{\varkappa}_{i,k} + \frac{1}{u+\gamma_{i,k}+ \hbar} \breve{\varkappa}'_{i,k}, \\
h_i(u) &\mapsto u^{-4m_1} \prod_{j=1}^{i-1} R_j(u-{\textstyle \frac{i-1-j}{2}\hbar}) \frac{\prod_{k=1}^{m_i} (u^2 - (\gamma_{i,k}+ \frac12 \hbar)^2)}{\prod_{l=1}^{m_{i-1}}(u-\gamma_{i-1,l})(u+\gamma_{i-1,l}+\hbar)}
\end{align}
defines a homomorphism $\breve{\Phi}^\lambda_\mu \colon \Tya \to \hdiffs$. 
\encor

\Rem
The symmetric formulation with square roots is somewhat more natural from the point of view of Gelfand--Tsetlin theory, see, e.g., \cite{GavKl} and \cite[Theorems 4.3, 6.1]{LP}. 
\enrem


\subsection{ABCD formulation}

Throughout this subsection, let $\mu = 0$. 
With a view to geometric applications, it is convenient to reformulate the twisted GKLO representations in terms of the ABCD presentation from \S \ref{sec: Sklyanin}. 

\Cor \label{cor: GKLO in ABCD form}
The twisted GKLO homomorphism $\Phi^\lambda_0 \colon \Tyazero \to {D}_{0,\hbar}^\lambda$ is uniquely determined by the following formulae: 
\begin{align*}
\widetilde{A}_i(u) &\mapsto u^{-2 m_i} \prod_{k=1}^{m_i} (u^2 - (\gamma_{i,k} + \half \h)^2), \\
\widetilde{B}_i(u) &\mapsto - u^{-2 m_i}  \sum_{k=1}^{m_i} \prod_{k\neq l=1}^{m_i} (u^2 - (\gamma_{i,l} + \half \h)^2)\Big( (u-\gamma_{i,k}-\half \h)\x_{i,k} + (u+\gamma_{i,k}+\half \h) \x'_{i,k} \Big), \\
\widetilde{C}_i(u) &\mapsto u^{-2 m_i} \sum_{k=1}^{m_i} \Big( \x_{i,k}(u-\gamma_{i,k}-\half \h) + \x'_{i,k}(u+\gamma_{i,k}+\half \h) \Big) \prod_{k\neq l=1}^{m_i} (u^2 - (\gamma_{i,l} + \half \h)^2). 
\end{align*} 
The symmetrized twisted GKLO homomorphism $\breve{\Phi}^\lambda_0 \colon \Tyazero \to \breve{D}_{0,\hbar}^\lambda$ is given by the same formulae, with replacements $\x_{i,k} \leftrightarrow \breve{\x}_{i,k}$ and $\x'_{i,k} \leftrightarrow \breve{\x}'_{i,k}$. 
\encor 

\Proof
This follows directly from Theorem \ref{thm: GKLO}, Corollary \ref{thm: GKLO sym}, as well as comparing formulae \eqref{eq: z in A terms} and \eqref{eq: ziu} (the former needs to be twisted by the action of the central element $c(u) \mapsto u^{-2\lambda_i}R_i(u)$ from \eqref{eq: central series}). The normalizing factor $u^{-2 m_i}$ is used to ensure $\widetilde{A}_i(u)$ is a series of the form $\widetilde{A}_i(u) = \sum_{k \geq 0} A_{i}^{(k)} u^{-k}$. 
\enproof 

We can now, at least partially, characterize the kernel of the twisted GKLO homomorphism. 

\Cor \label{cor: GKLO kernel}
The following elements are in the kernel of $\Phi^\lambda_0$: 
\begin{itemize}
    \item $\widetilde{A}_i^{(r)}, \widetilde{B}_i^{(r)}, \widetilde{C}_i^{(r)}$ for $r > 2m_i$. 
\end{itemize} 
\encor 

\Proof
The statement is immediate from Corollary \ref{cor: GKLO in ABCD form}. 
\enproof

\section{Geometric realisations}

Here we show that the shifted twisted Yangians from Section~\ref{sec-shifted twisted Yangians} and their truncations defined in Section~\ref{sec: tGKLO} quantise the Poisson structures discussed in Section~\ref{affgrassPoisson}.

\subsection{Quantisations via RTT generators}

It follows from the discussion in Section~\ref{sec-RTT} that 
$\dyazero/ \hbar \dyazero$
is generated as a $\mathbb{C}$-algebra by elements $s_{ij}^{(r)}$ with $1 \leq i,j \leq n$ and $r \geq 1$ subject to relations $s_{ij}^{(r)} = s_{ji}^{(r)}(-1)^r$. If $s_{ij}(z) := \sum_{r \geq 0} s_{ij}^{(r)} z^{-r}$ for $$s_{ij}^{(0)} := \begin{cases}
    1 & i=j \\
    0 & \text{otherwise}
\end{cases}$$ then the Poisson bracket on $\dyazero/ \hbar \dyazero$
induced as in Section~\ref{sec-setupcomm} is described by 
\begin{equation}\label{eq-twistedd RTT poisson relations}
\begin{aligned}
(u^2-v^2)\lbrace s_{ij}(u), s_{kl}(v) \rbrace =& (u+v)\left( s_{kj}(u) s_{il}(v) - s_{il}(u) s_{kj}(v) \right) \\
&- (u-v)\left( s_{ik}(u) s_{jl}(v) - s_{lj}(u) s_{ki}(v) \right)
\end{aligned}
\end{equation}
for variables $u,v$.
\begin{pro}\label{prop-twisted RTT isom}
    The map of $\mathbb{C}$-algebras 
    $$
     \dyazero/ \hbar \dyazero
     \xrightarrow{} \mathcal{O}(K_0 \backslash \operatorname{Gr}_0)
    $$
    given by $s_{ij}^{(r)} \mapsto \Delta_{ji}^{\tau,(r)}$ (recall Notation~\ref{notation-tau minors}) is a Poisson isomorphism.
\end{pro}
\begin{proof}
    Clearly the given map is an isomorphism of commutative rings. Recall from Section~\ref{affgrassPoisson} that $\lbrace-,-\rbrace$ denotes the Poisson bracket on $\mathcal{G}_0$, while $\lbrace - ,-\rbrace_\tau$ denotes the Poisson bracket on $\mathcal{O}(K_0 \backslash \mathcal{G}_0)$ obtained via $\Psi$ from the Dirac reduction of $\lbrace-,-\rbrace$ on $\mathcal{G}_0^{\tau=1}$. Lemma~\ref{lem-kwwypoissonstruc} and \eqref{eq-Dirac bracket} then combine to give
    \begin{equation}
\begin{aligned}
(u^2-v^2)\lbrace \Delta^\tau_{ij}(u), \Delta^\tau_{kl}(v) \rbrace_\tau =& (u^2-v^2)\lbrace \Delta_{ij}(u), \Delta_{kl}(v) \rbrace \circ \Psi \\
+& (u^2-v^2)\lbrace \Delta_{ji}(-u), \Delta_{kl}(v) \rbrace \circ \Psi\\
=&  (u+v)\left( \Delta^\tau_{il}(u) \Delta^\tau_{kj}(v)  
- \Delta^\tau_{kj}(u)\Delta^\tau_{il}(v) \right)  \\
-& (u-v) \left( \Delta^\tau_{jl}(-u) \Delta^\tau_{ki}(v) - \Delta^\tau_{ki}(-u)\Delta^\tau_{jl}(v) \right) \\
=&  (u+v)\left( \Delta^\tau_{il}(u) \Delta^\tau_{kj}(v)  
- \Delta^\tau_{kj}(u)\Delta^\tau_{il}(v) \right)  \\
-& (u-v) \left( \Delta^\tau_{lj}(u) \Delta^\tau_{ki}(v) - \Delta^\tau_{ik}(u)\Delta^\tau_{jl}(v) \right)
\end{aligned}
\end{equation}
Comparing with \eqref{eq-twistedd RTT poisson relations} shows that $s_{ij}^{(r)} \mapsto \Delta_{ij}^{(r)}$ as anti Poisson (i.e. maps the bracket on $\dyazero/ \hbar \dyazero$ onto $-1$ times the bracket on $\mathcal{O}(K_0\backslash \mathcal{G}_0)$). Since the association $s_{ij}^{(r)} \mapsto s_{ji}^{(r)}$ is an anti-automorphism of $\dyazero$, see \cite[Proposition 2.3.4]{Molev-book}, the claim follows.
\end{proof}

\begin{rem}
    The existence of the isomorphism in Proposition~\ref{lem-twisted identification} can also understood conceptually from the viewpoint of \cite{KWWY}. Specifically, \cite[Theorem 3.9]{KWWY} produces a Poisson isomorphism $\phi: \mathcal{O}(\operatorname{Gr}_0) \rightarrow \mathbf{Y}_0 / \hbar \mathbf{Y}_0$ where $\mathbf{Y}_0$ denotes the $\mathbb{C}[\hbar]$ form of the $\sln$ Yangian. This isomorphism identifies the universal matrix on $\operatorname{Gr}_0 \cong \mathcal{G}_0$ with the transpose of the matrix $T(z) \in \operatorname{Mat}(\mathbf{Y}_0 / \hbar \mathbf{Y}_0)[[z^{-1}]]$ of RTT generators. As a consequence, the isomorphism from Proposition~\ref{prop-twisted RTT isom} fits into the following commutative diagram
    $$
    \begin{tikzcd}
        \mathcal{O}(K_0 \backslash \operatorname{Gr}_0) \ar[r] \ar[d] & \ar[d]\mathcal{O}(\operatorname{Gr}_0) \\
        \dyazero/\hbar  \dyazero \ar[r]& \mathbf{Y}_0  / \hbar \mathbf{Y}_0
    \end{tikzcd}
    $$
    in which the lower horizontal embedding identifies  $ \dyazero/\hbar  \dyazero $ with the subalgebra of $\mathbf{Y}_0  / \hbar \mathbf{Y}_0$ generated over $\mathbb{C}$ by the coefficients of the entries inside $T(z)T^t(-z)$ \cite[Theorem 2.4.3]{Molev-book}.
\end{rem}
\subsection{Quantisations in the shifted setting}

\begin{pro}\label{lem-twisted identification}
	Let $\mu \in X_*(T)$ be dominant and recall the shift homomorphism ${}^{\mathsf{tw}}\iota(\mu): \dya \rightarrow \dyazero$. Then there is a commutative diagram of graded Poisson $\mathbb{C}$-algebras
	$$
	\begin{tikzcd}
		\mathcal{O}(K_0 \backslash \mathcal{G}_0 / \mathcal{U}^{-,\mu})  \ar[r] \ar[d,"\wr"]& \mathcal{O}(K_0 \backslash \mathcal{G}_0)  \ar[d,"\wr"] \\
		\dya / \h \dya   \ar[r,"{}^{\mathsf{tw}}\iota(\mu)"] & 	\dyazero / \h \dyazero 
	\end{tikzcd}
	$$
	whose right vertical arrow is the isomorphism in Proposition~\ref{lem-twisted identification}.
\end{pro}
\begin{proof}
    The image of $\dya / \hbar \dya$ inside $\dyazero / \hbar \dyazero$ under the shift homomorphism is Poisson generated by the images modulo $\hbar$ of the element $\hbar b_{i}^{(r)} \in \dyazero$ for $r > \langle \alpha_i^\vee,\mu \rangle$ and the $ \hbar z_{i}^{(r)} \in \dyazero$ for $r >0$ . On the other hand, Proposition~\ref{lem: Sk2} implies that the isomorphism in
    Proposition~\ref{lem-twisted identification} identifies the series
    \begin{equation}\label{eq-shifted generators}
    \hbar z_{i}(u) = \frac{A_{i-1}(u)A_{i+1}(u)}{A_i(u)^2}, \qquad  \hbar b_{i}(u) = B_i(u) A_i(u)
    \end{equation}
    for $A_i(u) = \sum_{r \geq 0} A_i^{(r)} u^{-r}, B_i(u) = \sum_{r >0} B_i^{(r)} u^{-r}$ and $A_i^{(r)},B_i^{(r)}$ as defined in Section~\ref{sec-ideal generators}.  It is easy to see that each coefficient in $A_i(u)$ lies inside the subring \linebreak $\mathcal{O}(K_0 \backslash \mathcal{G}_0 / 
    \mathcal{U}^{-,\mu})$ since the trailing principal minor of any $g \in G$ is invariant under left multiplication by $U^+$ and right multiplication by $U^-$. On the other hand, if $K_0x \in K_0 \backslash \mathcal{G}_0$ then $b_i(K_0x)$ is the $i,i+1$-th entry of the matrix $f \in \mathcal{U}_0^-$ obtained by factoring $\tau(x)x = edf$ as in \eqref{eq-Gauss factorisation}. It is easy to see that the coefficients in this series of degree $< \langle \alpha_i^\vee,\mu \rangle$ are invariant under the right action of $\mathcal{U}^{-,\mu}$. We conclude that, under the identification of Proposition~\ref{lem-twisted identification}, the image of $\dya / \hbar \dya$ lies inside $\mathcal{O}(K_0 \backslash \mathcal{G}_0 / \mathcal{U}^{-,\mu})$.
    
     It remains to show this inclusion is an equality. For this it suffices to observe that both  have Hilbert series for the loop grading given by
	 $$
	 \left( \prod_{\alpha^\vee \in \Delta^+} \prod_{q > \langle \alpha^\vee, \mu \rangle}^\infty \frac{1}{(1-t^q)} \right) \left(  \prod_{j=1}^\infty \frac{1}{(1-t^{2j})^{n-1}} \right). 
	 $$ 
For $\mathcal{O}(K_0 \backslash \mathcal{G}_0 / \mathcal{U}^{-,\mu})$ this is easily seen using the isomorphism
$$
K_0 \backslash \mathcal{G}_0 / \mathcal{U}^{-,\mu} \cong \left( \mathcal{U}_{-\mu}^+  \times \mathcal{T}_0 \times \mathcal{U}^-_\mu \right)^{\tau=1}
$$
in \eqref{eq-Psi}, while for $\dya/ \h \dya$, it follows from the description of the PBW basis in Proposition~\ref{pro: PBW}.
\end{proof}

\subsection{Quantisations of truncations}

Here we take $\mu =0$ and consider $\lambda \in X_*(T)^+$. Then Theorem~\ref{thm: GKLO} produces a surjection $\dyazero \rightarrow \dyatrzero$, and hence a surjection \linebreak $\dyazero / \hbar \dyazero \rightarrow \dyatrzero /\hbar \dyatrzero$.

\begin{thm}\label{thm-truncation}
    There is a commutative diagram
    $$
    \begin{tikzcd}
        \dyazero / \hbar \dyazero \ar[r,"\text{Prop~\ref{lem-twisted identification}}"] \ar[d]& \mathcal{O}(K_0 \backslash \mathcal{G}_0) \ar[d] \\
        \dyatrzero / \hbar \dyatrzero \ar[r] &  \mathcal{O}(\mathcal{S}_0^{\leq -w_0(2\lambda)}) 
    \end{tikzcd}
    $$
    whose bottom horizontal arrow is an isomorphism modulo the ideal of nilpotent elements in $\dyatrzero / \hbar \dyatrzero$.
\end{thm}
\begin{proof}
    Write $I$ for the Poisson ideal obtained as the kernel of the surjection $\mathcal{O}(K_0 \backslash \mathcal{G}_0) \cong \dyazero / \hbar \dyazero \rightarrow \dyatrzero /\hbar \dyatrzero$. From Corollary~\ref{cor: GKLO kernel} we see that 
    $$
    \Psi_0^\lambda(A_i^{(r)}) =0
    $$
    for $r > m_i := \langle \omega_i,2\lambda \rangle$.  Since $m_i = \langle \omega_i,-w_0(-w_0(2\lambda) \rangle$ it follows that $I$ contains the ideal described in Corollary~\ref{cor-ideal generators} when applied to the dominant coweight $-w_0(2\lambda)$. Thus, the vanishing locus $V(I) \subset K_0\backslash \mathcal{G}_0$ of $I$ is contained in $\mathcal{S}_0^{\leq -w_0(2\lambda)}$. If this containment were strict then $I$ could not contain any $A_i^{(m_i)}$. Indeed, Proposition~\ref{prop-nonvanishing} asserts that these functions are units in $\mathcal{O}(\mathcal{S}_0^{\leq -w_0(2\lambda)}$). However, it is clear from Corollary~\ref{cor: GKLO in ABCD form} 
    that $\Psi_0^\lambda(A_i^{(m_i)}) \neq 0$, so we are done.
\end{proof}

\begin{thm}
    Suppose that Conjecture~\ref{conj-reduced}. Then: 
    \begin{enumerate}
        
        \item The bottom horizontal arrow in Theorem~\ref{thm-truncation} is an isomorphism.
        \item $\dyatrzero$ is the quotient of $\dyazero$ by the two sided ideal generated by the $\hbar A_i^{(r)}$'s for $r> r_i = \langle \omega_i, 2\lambda \rangle$ and the $\hbar B_i^{(r_i+1)}$ for each $1 \leq i \leq n-1$. 
    \end{enumerate}
\end{thm}
\begin{proof}
    The argument is identical to that proving \cite[Theorem 4.10]{KWWY}. Let $K \subset \dyazero$ denote the two sided ideal described in the theorem. Then
    $$
    K/\hbar K \subset I \subset J_{0}^\lambda
    $$
    for $J_0^\lambda$ the ideal defining $\mathcal{S}_0^{\leq -w_0(2\lambda)}$ and $I$ the kernel of $\dyazero / \hbar \dyazero \rightarrow \mathcal{O}(\mathcal{S}_0^{\leq -w_0(2\lambda)})$. On the other hand, the proof of Theorem~\ref{thm-truncation} goes through with $I$ replaced by $K/\hbar K$ since Corollary~\ref{cor: GKLO kernel} clearly shows that 
    $$
    \Psi_0^\lambda(B_i^{(r)}) = 0
    $$
    for $r > \langle \omega_i,2\lambda \rangle$. Thus, $K/\hbar K$ has radical equal to $J_0^\lambda$. But Conjecture~\ref{conj-reduced} asserts that $K/\hbar K$ is already reduced, and hence $K/\hbar K = I =J_0^\lambda$ which gives part (1). Part (2) then follows from an application of Nakayama's lemma.
\end{proof}
\commentout{
\Pro
We have: 
\[ 
\mathcal{L}_\mu^{\leq \lambda} \subseteq V(I_\mu^\lambda) \subseteq V(J_\mu^\lambda). 
\]
\enpro 

\Proof
The second inclusion follows directly from Lemma \ref{lem: tdelta under Psi} and the fact that both $J_\mu^\lambda$ and $I_\mu^\lambda$ are Poisson ideals. 

The vanishing locus of $I_\mu^\lambda$ is a union of symplectic leaves. Assuming (*), if the first inclusion does not hold, then 
\[
V(I_\mu^\lambda) \cap \mathcal{L}_\mu^{\leq \lambda} = \bigcup_j \mathcal{L}_\mu^{\leq \nu_j} =:X, 
\]
for some $\nu_j < \lambda$. 
For each $j$, there exists $i$ such that $\langle \nu_j, \omega_{i^*} \rangle < \langle \lambda, \omega_{i^*} \rangle = m_i$. By \cite[Proposition 2.4]{KWWY}, $M = \prod_i \tdelta_{i,i}^{(2m_i)}$ vanishes on $X$. Hence, for some $k>0$, $M^k \in I^\lambda_\mu$. However, under the twisted GKLO homomorphism, $\Phi^\lambda_\mu(\tdelta_{i,i}^{(2m_i)}) = \gamma_{i,1}^2 \cdots \gamma_{i,m_i}^2$. Hence $M^k$ is mapped to a monomial in the $\gamma_{i,k}$, which is non-zero, and so $M^k \notin I^\lambda_\mu$, resulting in a contradiction. 
\enproof
} 



\newcommand{\etalchar}[1]{$^{#1}$}
\providecommand{\bysame}{\leavevmode\hbox to3em{\hrulefill}\thinspace}
\providecommand{\MR}{\relax\ifhmode\unskip\space\fi MR }
\providecommand{\MRhref}[2]{%
  \href{http://www.ams.org/mathscinet-getitem?mr=#1}{#2}
}
\providecommand{\href}[2]{#2}

\end{document}